# A minimization principle behind the diffusion bridge of diurnal fish migration

Hidekazu Yoshioka[1,*]

[1] Graduate School of Advanced Science and Technology, Japan Advanced Institute of Science and Technology, 1-1 Asahidai, Nomi, Ishikawa, Japan

[*] Corresponding author: yoshih@jaist.ac.jp, ORCID: 0000-0002-5293-3246

**Abstract**

Fish migration is a mass movement that affects the hydrosphere and ecosystems. While it occurs on multiple temporal scales, including daily and intraday fluctuations, the latter remains less studied. In this study, for a stochastic differential equation model of the intraday unit-time fish count at a fixed observation point, we demonstrate that the model can be derived from a minimization problem in the form of a stochastic control problem. The control problem assumes the form of the Schrödinger Bridge but differs from classical formulations by involving a degenerate diffusion process and an objective function with a novel time-dependent weight coefficient. The well-posedness of the control problem and its solution are discussed in detail by using a penalized formulation. The proposed theory is applied to juvenile upstream migration events of the diadromous fish species *Plecoglossus altivelis altivelis* commonly called *Ayu* in Japan. We also conduct sensitivity analysis of the models identified from real data.

**Keywords**

Diurnal fish migration, Schrödinger Bridge, Optimal control, Penalty method, *Plecoglossus altivelis altivelis*

***Statements & Declarations***

**Acknowledgments:** The author would like to express his gratitude to Japan Water Agency, Ibi River, and Nagara River General Management Office for providing the valuable fish migration data. The author is very grateful for the helpful comments from the reviewers.

**Funding:** This work was supported by the Japan Science and Technology Agency (PRESTO No. JPMJPR24KE).

**Competing Interests:** The authors have no relevant financial or nonfinancial interests to disclose.

**Data Availability:** The data will be made available upon reasonable request from the corresponding author.

**Declaration of Generative AI in Scientific Writing:** The authors did not use generative AI for the scientific writing of this manuscript.

**Contribution:** All parts of this manuscript were prepared by the sole author.

## 1. Introduction

In this paper, we address the theoretical reasoning of a singular coefficient in a mathematical model for describing fish migration observed at a fixed point. We address this issue from a viewpoint of a novel stochastic control approach and apply the model to real data.

### 1.1 Problem background

Migration is a biological phenomenon involving synchronous movement of animals. Migratory fish species are model organisms for studying migration because they cover various migration forms [1]. Furthermore, migratory fish species play indispensable roles in circulating aquatic environments and ecosystems on Earth. They transport nutrients between different water bodies [2] and understanding the impact of their role under climate change is important [3]. The degradation of the longitudinal connectivity of river caused by physical barriers hinders their movements to suitable habitats [4] and can disrupt their life cycles [5]. Migratory fishes are also economically important as fishery resources, necessitating close monitoring of their population dynamics under climate change and anthropogenic pressure [6,7]. Therefore, the study of fish migration is a key research topic in many aspects.

Fish migration involves a wide range of timescales, from seasonal to intraday. A common timescale in the literature on fish migration, particularly for those studying fish migration along rivers, is daily to monthly, and migrating fish are typically counted at fixed observation point(s). Such a survey can be effective for investigations of the seasonal migration behavior of the target fish species because tracking changes of daily fish count enables assessing annual fishway evaluation [8], differentiating migration strategies among distinct fish species [9], analyzing influences of physical barriers [10], detecting changes of life cycle events, such as spawning [11], recruitment [12], and smolting [13].

By contrast, investigating intraday migration dynamics is essential for understanding the biological and ecological aspects of migratory fishes at finer temporal scales, because daily count data lack the resolution to capture hourly or sub-hourly patterns. The sub-daily run size of the post-larval migration of a goby was studied using a fixed net system [14]. Acoustic telemetry has been employed to study the influence of intraday tidal gate operation on the passage efficiency of sea lamprey [15] and the migration of Atlantic salmon through lakes [16]. Imaging sonar has been used to differentiate the intraday frequency of migration events among fish species [17]. Radio telemetry has been used to evaluate the passage efficiency of Atlantic salmon through fishways [18]. The hourly relationship between the spawning activity of alewife and herring and predator presence has been statistically assessed using an automated acoustic detection method [19]. High-resolution imaging sonar has been employed to clarify the daily and intraday migration dynamics of herring in coastal rivers [20]. Hourly fish count data were obtained using a video recording system, which suggested that the unit-time fish count fluctuates significantly with size and frequency, depending on the direction of migration [21]. The 30-min fish count data of the spawning population of Asian salmonids were estimated using image sonar [22]. The 10-min fish count of *Plecoglossus altivelis altivelis*, commonly called *Ayu* in Japan, during its upstream migration was monitored by a video system to investigate the fish count as a function of discharge and salinity [23].

Despite ongoing efforts, intraday fish migration is still less theoretically understood compared to seasonal and daily fish migration patterns. One plausible reason is that manually and/or visually tracking migrating individual fish, which are fast moving objects, is a challenging task that is labor and cost intensive. Recently, Yoshioka [24] proposed a stochastic differential equation (SDE) model to describe the diurnal migration phenomenon of *Ayu* based on a bridge, a stochastic process with pinned initial and terminal conditions. In the context of the diurnal fish migration, the state variable of the model is the unit-time fish count at a fixed observation point in the river, and the initial and terminal times correspond to sunrise and sunset, respectively. The model has been applied to real 10-min migration data and can reveal the fluctuation and intermittency of fish count [24]; however, the model was derived only heuristically; hence, the origin of its functional form is still unclear. To improve the applicability of this model and explore advanced models, its underlying mathematical laws and biological mechanisms should be clarified. With this contribution, a deeper theoretical analysis of the model would become feasible and potentially provide a unified viewpoint that also applies to other migratory fish species.

### 1.2 Aim and contribution

This study aims to provide a novel theoretical background of the diffusion bridge of Yoshioka [24] for heuristically modeling intraday fish counts. The contributions of this study are as follows:

We explicitly relate this diffusion bridge to a minimization problem of an energy functional, identifying the bridge's coefficient as the optimal control variable. This connection clarifies the functional form and properties of the bridge, suggesting that the minimization principle underlies intraday fish migration. Therefore, the main tasks of this study are to find an energy functional, called objective function in control theory, whose minimization yields a diffusion bridge and to investigate its properties. While earlier studies have suggested that animal migration may follow some optimality principles—balancing risks and rewards [25-27]—a clear mathematical formulation has not been established for intraday fish migration. We address this gap by a mathematical approach.

We adopt stochastic control formalism [28] well-suited for modeling and analyzing stochastic differential equations (SDEs). The proposed control formalism is based on Yoshioka [29] who discussed a problem with a limited functional form of the objective function and justified diffusion bridges only for one model case; it was therefore not practical and was not able to cover the whole models discussed in Yoshioka [24]. To overcome this critical issue, we extend the objective function so that a boarder range of diffusion bridges are covered. Other existing models of animal migration in a related direction are the mean-field games of diel vertical migration of single [30-32] and multiple aquatic species [33]. These and our studies share the common principle that biological movement emerges from an optimality principle, while the mathematical difference is that the existing studies discussed SDEs with a constant diffusion coefficient, whereas ours uses a non-Lipschitz continuous and degenerate diffusion coefficient. Degenerate diffusion, which poses a technical difficulty, is essential for a diffusion bridge to physically guarantee the non-negativity of fish count. Optimal control and dynamic game models for determining migration timing have been discussed [34-36]; however, in principle they cannot address temporally fine migration events.

Our control problem resembles the Schrödinger Bridge, which is a control problem with a prescribed terminal condition; therefore, it is not a standard stochastic control problem but harmonizes with modeling bridges in various fields, including biology, physics, and machine learning [37-42]. We resolve the technical issue of the terminal condition by formulating a penalized problem in which the terminal condition is softly constrained in an objective function to be optimized [e.g., 43,44]. Our penalized problem admits a closed-form solution because it exploits the affine nature [45] of the diffusion bridge, and the Hamilton–Jacobi–Bellman (HJB) equation associated with the problem admits an explicit smooth solution. A verification argument for smooth solutions (e.g., Theorem 3.5.2 in Pham [28]) proves that the solution is indeed the value function of the penalized problem. Considering a strong penalization limit yields the value function and optimal control of the original control problem under certain conditions. We study these conditions in detail and provide insights into fish migration phenomena. We emphasize that the objective function in our control formalism involves a time-dependent weight coefficient being different from Yoshioka [29], with which we can cover realistic diffusion bridges in applications. In particular, we derive a sharp condition under which the optimally controlled dynamics equal the diffusion bridge of fish migration.

Finally, we apply the minimization principle to the intraday fish count data of spring upstream migration events of juvenile *Ayu* in Japan and discuss the minimization principle behind their migration. Models with time-dependent sources and volatility are identified from the real 10-min fish count data of *Ayu* focusing on a case where the time dependence of the weight in the objective function for the minimization principle plays a role. We also conduct sensitivity analysis of the identified models. Consequently, this study provides a new theoretical foundation for a model of intraday fish migration, along with an application.

Note that high-resolution (e.g., 10-min) data of intraday fish count are rare by themselves due to technical reasons for implementing an appropriate observation system. The author could fortunately obtain the data used in this paper in 2025 thanks to the kindness of Japan Water Agency, Ibi River, and Nagara River General Management Office; the Nagara River is a valuable example where high-resolution fish count data have been collected. This is a possible reason mathematical models, including SDEs, have not been investigated for temporally-fine fish count data. We therefore consider that the model in Yoshioka [24] is a founding model for describing high-resolution unit-time fish count at a fixed location in a river. This paper contributes to deeper understandings about this model.

The remainder of this paper is organized as follows. **Section 2** reviews the diffusion bridges for intraday fish migration. **Section 3** presents a control problem whose solution potentially gives a diffusion bridge. **Section 4** applies the proposed minimization principle to *Ayu* in Japan along with sensitivity analysis of identified models. **Section 5** concludes the study and presents its perspectives. **Appendix** presents proofs of propositions and auxiliary results.

## 2. Diffusion bridge

We review the diffusion bridge and discuss its connection to stochastic control.

### 2.1 The model

In this study, time $t$ is a non-negative parameter, and the unit-time fish count at a fixed observation point, approximated as a continuous variable, is denoted as $X_t$ at time $t$. Physically, the unit-time fish count must be non-negative. We focus on the diurnal migration phenomenon of fish, such that the fish count becomes positive only during sunrise (time 0) and sunset (time $T > 0$) in a day. The 1-D standard Brownian motion is denoted $B_t$ at time $t$. In the sequel, when there will be no confusion, we omit subscripts with respect to time when we assume a stochastic process but not its value, for example $B$ but not $(B_t)_{0 \le t \le T}$

Yoshioka [24] proposed the Itô's SDE to describe the diurnal fish count dynamics:

$$\mathrm{d}X_t = (a - h_t X_t)\mathrm{d}t + \sigma\sqrt{h_t X_t}\mathrm{d}B_t, \quad 0 < t < T \tag{1}$$

subject to homogeneous initial and terminal conditions $X_0 = X_T = 0$. Here, $a > 0$ is the source coefficient representing the supply of migrants, $\sigma > 0$ is the volatility representing the intensity of fluctuation in the fish count, and $h_t$ is the reversion rate at time $t$, which is a positive function of time. The constraints $X_0 = X_T = 0$ represent the assumption that fish migration is diurnal and occurs only during the daytime. We consider that the supply of migrants from the downstream river reach contributes to the increase of the unit-time fish count deterministically as described in the SDE (1), and the contribution is assumed to be constant here for modeling simplicity. Further detailed clarifications of the source $a$ (e.g., when and how migrants are supplied from downstream) need a mechanistic explanation, which cannot be addressed by our conceptual model, and we consider that this point is a limitation of its explainability.

The appearance of the coefficient $h$ in both the drift (first) and diffusion (second) terms in the SDE (1) is based on the ansatz that fish migration is modulated by a biological clock, where the modified time is expressed by the primitive $H_t$ of $h_t$, i.e.,

$$-h_t X_t \mathrm{d}t + \sigma\sqrt{h_t X_t}\mathrm{d}B_t = -X_t \mathrm{d}(H_t) + \sigma\sqrt{X_t}\mathrm{d}(B_{H_t}), \tag{2}$$

implying that the diurnal behavior controlled by $h$ is mathematically a time change of the SDE. Mathematically, this is an equality to be understood in the sense of law. Biologically, this equation is an assumption that both reversion and fluctuation are influenced by some biological clock while the source is not. More specifically, we assume that the source (supply of migrants from downstream) is an external input to the fish count at the observation point that is the reason that the term including $a$ is not considered in (2).

### 2.2 Model properties

Yoshioka [24] heuristically proposed the SDE (1) as a minimal model that complies with the following requirements for describing diurnal fish migration:

***Assumptions***

(a) The SDE (1) admits a unique pathwise solution that is non-negative for $0 < t < T$ under suitable conditions irrespective to values of $a, \sigma > 0$.

(b) The initial and terminal conditions $X_0 = X_T = 0$ are satisfied with probability 1 under suitable conditions irrespective to values of $a, \sigma > 0$.

(c) Parameters $a, \sigma$ and coefficient $h$ can be identified from real data sets. In the simplest case, the coefficient $h$ has only one parameter to be fitted.

(d) The moments and cumulants of $X$ are found in closed forms because they are affine processes.

(e) The SDE (1) can be numerically computed using a modern numerical method that preserves the non-negativity of numerical solutions.

The first (a) and second items (b) are for the well-posedness, and they hold true with $h$, which is bounded and continuously differentiable in $(0,T)$ such that

$$\frac{c}{T-t} \leq h_t \leq \frac{c}{T-t} + c', \quad 0 < t < T \tag{3}$$

with some constants $c, c' > 0$, i.e., $h$ that blows up at the rate $(T-t)^{-1}$ near $T$. This unboundedness property of $h$ is the key to well-pose the model; the condition (3) means that $h$ should be neither too large nor too small for the well-posedness of the model. The third item (c) is based on practicality and satisfied by specific $h$, including the simplest case $h_t = \dfrac{c}{T-t}$ with $c > 0$. The fourth item (d) is valuable both theoretically and practically, as it eliminates the need to simulate sample paths when the focus is on moments or cumulants. The final item (e), based on the dynamics-preserving discretization method [e.g., 46], is important in applications because it allows generating sample paths and thus, computing path-wise quantities that cannot be evaluated using moments.

The key of the SDE (1), in addition to the coefficient $h$, is its non-Lipschitz continuous and degenerate diffusion coefficient, which is proportional to the square root $\sqrt{X_t}$. This square root nature of the diffusion balances the well-posedness and model usability due to the following reasons. First, the coefficient $\sqrt{X_t}$ guarantees the non-negativity of the unique solution to the SDE (1); intuitively, with this diffusion coefficient, sample paths of $X$ possibly touch but do not cross the value 0. Second, the infinitesimal generator associated with $X$ has affine coefficients due to the square-root diffusion. This property is exploited in the control formalism presented in this paper. Third, both low- and high-volatility regimes can be simulated by specifying sufficiently small and large volatility values, respectively. Here, the low-volatility regime indicates that sample paths of $X$ do not reach 0. By contrast, the high-volatility regime indicates that sample paths of $X$ reaches 0 within a finite time with probability 1. The latter occurs if the following condition is satisfied [24]:

$$a < \frac{\sigma^2}{2} \inf_{t \in (0,T)} h_t . \tag{4}$$

See **Figure 1** for samples paths with and without Feller condition.

The condition (4) is a version of the (violation of) Feller condition, which determines whether the target process hits a boundary. The case in which a Feller condition is violated tends to be avoided in literature possibly because of the theoretical complexity that is not encountered when the Feller condition is satisfied. Nevertheless, the violation becomes essential in certain financial applications [47] and theoretical studies [48,49]. In our context, (4) indicates that fish migration is essentially in a high-volatility regime, and hence, the Feller condition is violated; see also **Section 4**.

In summary, the square-root diffusion coefficient is a reasonable choice in the sense that it guarantees the well-posedness of the model, preserves the non-negativity of its solution, complies with the biological background (2), is able to generate intermittent sample paths as in the real data [24] with an explicit threshold (Feller condition), and is the one whose associated HJB equation is analytically solvable as demonstrated later in this paper. For example, a modeling alternative is a diffusion coefficient proportional to the solution, but the intermittency condition would become different, even if it exists. With a more nonlinear and complex diffusion coefficient, it is unclear whether the well-posedness and non-negativity of the model are satisfied.

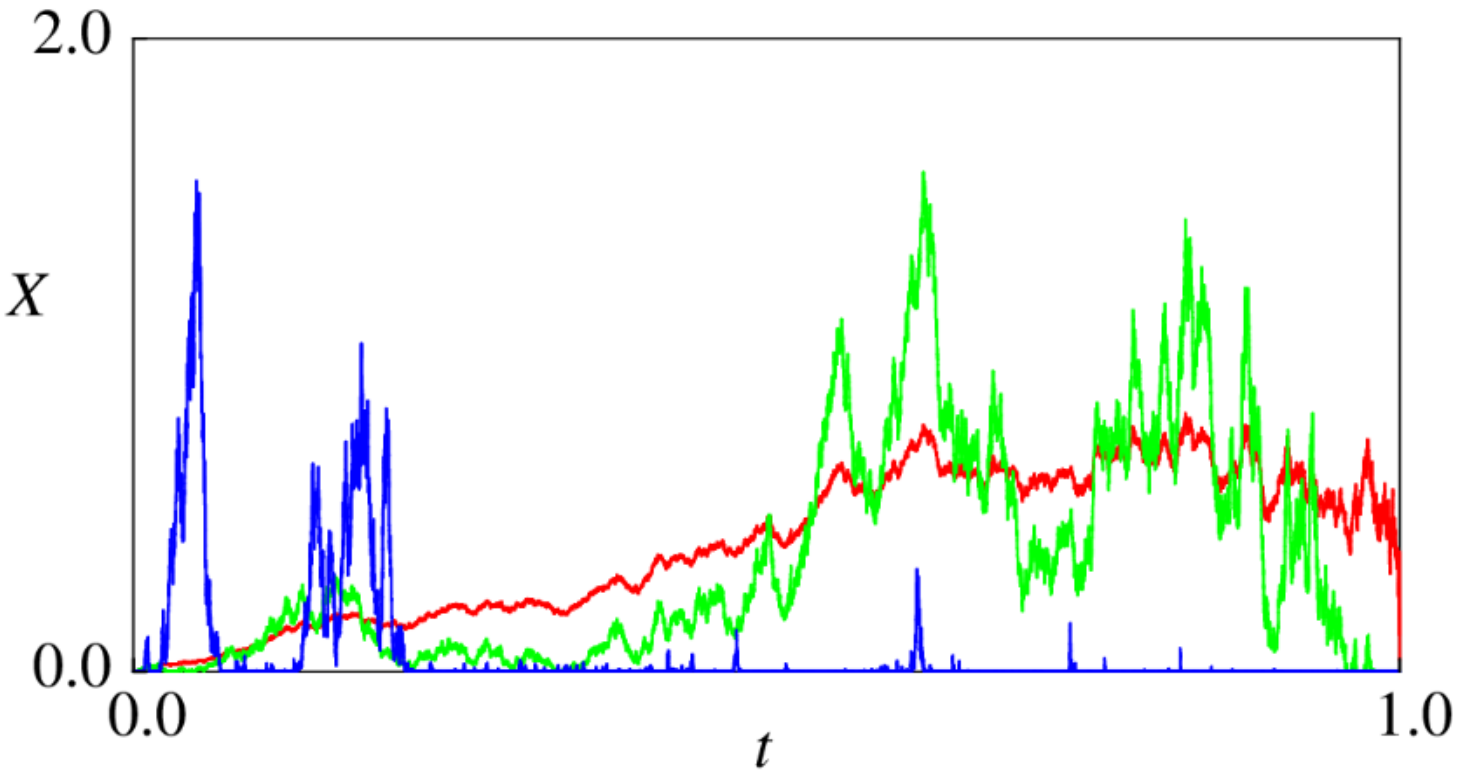

**Figure 1.** Sample paths of the SDE (1) for low- and high-volatility regimes. The terminal time is 1. We choose $a=1$ and $h_t = 0.5(1-t)^{-1}$: $\sigma = 0.5$ (red), $\sigma = 2.0$ (green), and $\sigma = 8.0$ (blue).

### 2.3 Control perspective

As discussed in the previous section, the SDE (1) depends on the coefficient $h$, particularly its blow-up speed $(T-t)^{-1}$ near the time $T$. This blow-up speed has been specified to well-pose the SDE (1) along with the constraints $X_0 = X_T = 0$. Here, a natural question arises: how can we derive the coefficient $h$ that satisfies these conditions from an optimality principle without adding assumptions? We address this by formulating a control problem where the coefficient $h$ is treated as the optimized control variable.

## 3. Control problem

### 3.1 Schrödinger Bridge

We explore the following control problem (**Problem O**): *find the minimum value of the objective function* $J$ *, namely*

$$\Phi(t,x)=\inf_{u\in U}J(t,x;u)=\inf_{u\in U}\mathbb{E}\left[\int_t^T\frac{1}{m+1}\frac{1}{w_s}X_s^{(u)}u_s^{m+1}\mathrm{d}s\,\middle|\,X_t^{(u)}=x\right],\quad 0<t<T,\quad x\geq 0 \tag{5}$$

*subject to the controlled SDE*

$$\mathrm{d}X_t^{(u)}=\left(a-u_tX_t^{(u)}\right)\mathrm{d}t+\sigma\sqrt{u_tX_t^{(u)}}\mathrm{d}B_t,\quad 0<t<T \tag{6}$$

*with the initial and terminal constraints*

$$X_0^{(u)}=X_T^{(u)}=0. \tag{7}$$

Here, $m>0$ is a power constant, $(w_t)_{0\leq t<T}$ is a positive and continuously differentiable function in $[0,T)$ with $w_0=1$, and $U$ is the admissible set of control variables $u$ given by

$$U=\left\{(u_t)_{0\leq t<T}\,\middle|\,u\text{ is measurable and adapted to }\mathbb{F}\text{ , }u_t\geq 0,\text{ and }J(\cdot,\cdot,u)<+\infty\right\}, \tag{8}$$

where $\mathbb{F}$ is a natural filtration generated by $(B_t)_{0\leq t\leq T}$. The set $U$ is not empty because the null control $u\equiv 0$ belongs to it.

**Problem O** can be equivalently rewritten with a free endpoint as follows:

$$\Phi(t,x)\equiv\inf_{u\in U}\mathbb{E}\left[\int_t^T\frac{1}{m+1}\frac{1}{w_s}X_s^{(u)}u_s^{m+1}\mathrm{d}s+\mathbb{I}\left(X_T^{(u)}>0\right)X_T^{(u)}\,\middle|\,X_t^{(u)}=x\right],\quad 0<t<T,\quad x\geq 0 \tag{9}$$

subject to the SDE (6) but now without constraints (7). Here, $\mathbb{I}(y>0)$ ($y\geq 0$) is an indicator function such that $\mathbb{I}(y>0)=+\infty$ if $y>0$ and $\mathbb{I}(y>0)=0$. Indeed, controls that violate the terminal condition $X_T^{(u)}=0$ are never optimal.

Our control problem is conceptual, that is, it is not fully mechanistic; nevertheless, it potentially accounts for some ecological aspects of fish migration as discussed below. The decision-maker in **Problem O** is the migrating fish population as an aggregation of individual fish, assuming that there are social cues that trigger and modulate their migration [50-52]. The controlled SDE (6) has the same form as the bridge (1) except that the reversion rate $h$ is replaced by a control $u$. The objective function $J$, which is minimized by choosing controls, measures the net unit time cost of migration $\frac{1}{m+1}\frac{1}{w_s}X_s^{(u)}u_s^{m+1}$ at each time $s$ conditioned on the event $X_t^{(u)}=x$. This unit-time cost implies that the cost is weighted by the coefficient $w_s$, proportional to the unit-time fish count $X_s^{(u)}$, and is a convex and increasing function of the control $u_s$. The parameter $m$ modulates sensitivity of the cost about the control $u_s$, in a way that a larger $m$ more strongly suppresses rapid variations of the unit-time fish count ($\frac{1}{m+1}$ is simply a normalization factor). Larger values of $w_s$ imply smaller costs of migration, which would implicitly depend on environmental conditions, such as water temperature, water depth, flow velocity, light level, and photoperiod, which are not controllable by fish, and the cost of migration would decrease if their ranges are suitable for migrating fish [53-57]. Among them, hydrodynamic costs may be a dominant factor that influences the physiological cost of fish migration [58-60], potentially affected by water temperature [61] and school shape [62,63]. In this view, it may be natural and theoretically simpler to select the parameter value $m = 1$ such that the objective functional resembles kinetic energy. However, choosing $m = 1$ is not necessarily permissible considering the well-posedness of the minimization problem (see **Section 3.3**). Thus, considering the generic $m > 0$ is essential for the modeling in this study. Finally, the earlier control formalism of Yoshioka [29] corresponds to the constant $w$ case, which has a quite limited applicability as discussed in **Remark 2** presented later.

In (9), proportionality with respect to $X_s^{(u)}$ implies that the migration of a larger population is more costly than that of a smaller population because of the increasing risk of predation. For example, this corresponds to the situation where juveniles of migratory fish can be more susceptible to predation by piscivorous birds and terrestrial predators [64-67]. The proposed objective function is assumed to inclusively describe the physical, biological, and ecological factors discussed above.

Mathematically, **Problem O** is not a standard stochastic control problem but a Schrödinger Bridge because both the initial condition $X_0$ and the terminal condition $X_T$ are prescribed; the latter is represented by a singular terminal cost in (9). This prevents us from applying existing stochastic control theory [e.g., 28] to our control problem. To resolve this theoretical difficulty, we consider the following penalized problem, where the terminal constraint is softly incorporated into the objective function, as employed in modeling bridges [42,43] (**Problem P**): *find the minimum value of the penalized objective function* $J_\eta$*, namely*

$$\Phi_\eta(t,x) \equiv \inf_{u\in U} J_\eta(t,x;u) = \inf_{u\in U} \mathbb{E}\left[\int_t^T \frac{1}{m+1}\frac{1}{w_s} X_{\eta,s}^{(u)} u_s^{m+1} \mathrm{d}s + \frac{1}{\eta} X_{\eta,T}^{(u)} \middle| X_{\eta,t}^{(u)} = x\right],\quad 0 < t < T,\quad x \geq 0, \quad (10)$$

where $X_\eta^{(u)}$ satisfies the SDE (6) with the initial condition $X_{\eta,0}^{(u)} = 0$. Here, $\eta > 0$ is the penalty parameter and the added term $\frac{1}{\eta} X_{\eta,T}^{(u)}$ in (10) measures the deviation between the desired constraint 0 and the attained value $X_{\eta,T}^{(u)}$. We have the following pointwise convergence:

$$\frac{1}{\eta} y \underset{\eta\to+0}{\to} \mathbb{I}(y > 0) \quad \text{at each} \quad y \geq 0. \quad (11)$$

Therefore, **Problem P** is a relaxed version of **Problem O**, and we expect that the former gets closer to the latter as penalization becomes stronger, i.e., as $\eta \to +0$. In contrast to the literature that employs a quadratic penalty term proportional to $\left(X_{\eta,T}^{(u)}\right)^2$ [43], we use the affine penalty $X_{\eta,T}^{(u)}$ because the controlled process $X_{\eta}^{(u)}$ remains non-negative with probability 1. Moreover, the affine penalty term harmonizes with the affine nature of the controlled SDE.

***Remark 1*** A significant difference between the classical Schrödinger Bridges and our formulation is that the former typically involves additive noise in the state dynamics [e.g., 37-42], whereas ours is driven by a multiplicative noise. Another difference is that our method employs non-quadratic, state-dependent (proportional to $X_{\eta,s}^{(u)}$), and time-dependent weights ( $w_s$ ) in the objective function. Using a time-dependent weight $w_s$ is crucial for obtaining a sufficiently wide applicability range for the proposed control framework; see also **Remark 2**.

### 3.2 Penalized problem

**Problem P** assumes the form of a control problem to which the dynamic programming principle and verification theorem apply; hence, it is solvable within the classical framework of optimal control. Specifically, we use a verification argument which states that if the following HJB equation admits a classical (sufficiently smooth) solution, then it is the value function $\Phi_{\eta}$ of **Problem P**:

$$\frac{\partial \phi(t,x)}{\partial t} + a\frac{\partial \phi(t,x)}{\partial x} + x\inf_{v\ge 0}\left\{ -v\frac{\partial \phi(t,x)}{\partial x} + \frac{\sigma^2}{2}v\frac{\partial^2 \phi(t,x)}{\partial x^2} + \frac{1}{w_t}\frac{v^{m+1}}{m+1} \right\} = 0\,,\quad 0<t<T\,,\quad x\ge 0 \tag{12}$$

subject to the terminal condition $\phi(T,x) = \frac{1}{\eta}x$ ( $x \ge 0$ ). Moreover, the optimal control $u = u_{\eta}^{*}$ as the minimizer in (10) is given by

$$u_{\eta,t}^{*} = \arg\min_{v\ge 0}\left\{ -v\frac{\partial \Phi_{\eta}\left(t,X_{\eta,t}^{*}\right)}{\partial x} + \frac{\sigma^2}{2}v\frac{\partial^2 \Phi_{\eta}\left(t,X_{\eta,t}^{*}\right)}{\partial x^2} + w_t\frac{v^{m+1}}{m+1} \right\},\quad 0<t<T\,, \tag{13}$$

where $X_{\eta,t}^{*}$ means $X_{\eta,t}^{\left(u_{\eta}^{*}\right)}$. Therefore, the main task is to solve the HJB equation (12).

The following **Proposition 1** explicitly provides a solution to (12) and the corresponding optimal control (13). The condition (14) gives an design criterion to well-pose the control problems.

***Proposition 1 (Optimality result of Problem P)***

*Assume that*

$$0 \le w_t^{\frac{1}{m}} \le c''(T-t)^{-\alpha} \quad \textit{for any} \quad 0\le t<T \quad \left(\textit{this implies } \int_0^T w_s^{\frac{1}{m}}\mathrm{d}s < +\infty\right) \tag{14}$$

*with some constants* $c''>0$ *and* $\alpha \in [0,1)$. *Then, it follows that*

$$\Phi_{\eta}(t,x) = A_{\eta,t}x + C_{\eta,t}\,,\quad 0<t<T\,,\quad x\ge 0\,, \tag{15}$$

*where*

$$A_{\eta,t} = \left( \eta^{\frac{1}{m}} + \frac{1}{m+1} \int_t^T w_s^{\frac{1}{m}} \mathrm{d}s \right)^{-m}, \quad 0 < t < T \tag{16}$$

*and*

$$C_{\eta,t} = a \int_t^T A_{\eta,s} \mathrm{d}s, \quad 0 < t < T. \tag{17}$$

*Moreover, the corresponding optimal control is given by*

$$u_{\eta,t}^* = w_t^{\frac{1}{m}} \left( A_{\eta,t} \right)^{\frac{1}{m}}, \quad 0 < t < T. \tag{18}$$

### 3.3 Limit to the original control problem

By obtaining the solution to the penalized problem (**Proposition 1**), the next task is to analyze its limit $\eta \to +0$ and check whether the limit is non-trivial (i.e., $\Phi_\eta$ is not identical to $+\infty$ ). Another important task is to determine the condition under which $u^* = h$, where $u^*$ is the minimizer of **Problem O**, so that the control problem is linked to the diffusion bridge. To investigate these issues, we must specify the form and/or behavior of the coefficient $w$.

By (16) and (17), we have

$$A_t \equiv \lim_{\eta \to +0} A_{\eta,t} = \left( \frac{1}{m+1} \int_t^T w_s^{\frac{1}{m}} \mathrm{d}s \right)^{-m}, \quad 0 \le t < T \tag{19}$$

and

$$C_t \equiv \lim_{\eta \to +0} C_{\eta,t} = a \int_t^T A_s \mathrm{d}s = a(m+1)^m \int_t^T \left( \int_s^T w_q^{\frac{1}{m}} \mathrm{d}q \right)^{-m} \mathrm{d}s, \quad 0 \le t < T. \tag{20}$$

Considering (19) and (20), we infer that the optimal control of **Problem O** is given by

$$u_t^* = w_t^{\frac{1}{m}} \left( A_t \right)^{\frac{1}{m}} = (m+1) w_t^{\frac{1}{m}} \left( \int_t^T w_s^{\frac{1}{m}} \mathrm{d}s \right)^{-1}, \quad 0 < t < T. \tag{21}$$

If (21) holds true, then we can determine the coefficient $w$ by equating $u^*$ and $h$. If so, assuming that $w$ and $h$ are bounded and continuously differentiable in $(0,T)$, by $u^* = h$ and (21), we obtain

$$h_t \int_t^T w_s^{\frac{1}{m}} \mathrm{d}s = (m+1) w_t^{\frac{1}{m}}, \quad 0 < t < T. \tag{22}$$

Then, we derive the linear ordinary differential equation to determine $z_t = w_t^{\frac{1}{m}}$ when $h$ is given (differentiating both sides of (22) using $t$):

$$\frac{\mathrm{d}z_t}{\mathrm{d}t} = \left( \frac{1}{h_t} \frac{\mathrm{d}h_t}{\mathrm{d}t} - \frac{1}{m+1} h_t \right) z_t, \quad 0 < t < T, \tag{23}$$

which is solvable explicitly along with $w_0 = 1$:

$$w_t = z_t^m = \left( \frac{h_t}{h_0} \right)^m e^{-\frac{m}{m+1} \int_0^t h_s \mathrm{d}s}, \quad 0 \le t < T. \tag{24}$$

The formula (24) is the key to deriving the diffusion bridge (1) from the control problem and suggests how to design the $w$ to obtain the desired reversion rate $h$. For example, if $h_t = \dfrac{c}{T-t}$ with $c > 0$, then we obtain

$$w_t = \left(1 - \frac{t}{T}\right)^{\frac{m}{m+1}(c-(m+1))}, \quad 0 \le t < T, \tag{25}$$

and then the assumption (14) is satisfied since

$$\int_0^T w_s^{\frac{1}{m}} \mathrm{d}s = \int_0^T \left(1 - \frac{s}{T}\right)^{\frac{c}{m+1}-1} \mathrm{d}s = T\frac{m+1}{c} < +\infty . \tag{26}$$

To have the boundedness of $\Phi$, we need $C_0 < +\infty$. By (20), an elementary calculation yields

$$C_0 = a(m+1)^m \int_0^T \left( \int_s^T \left(1 - \frac{q}{T}\right)^{\frac{c}{m+1}-1} \mathrm{d}q \right)^{-m} \mathrm{d}s = \frac{ac^m}{T^m} \int_0^T \left(1 - \frac{s}{T}\right)^{\frac{-m}{m+1}c} \mathrm{d}s = \begin{cases} +\infty & \left(c \ge 1 + \dfrac{1}{m}\right) \\ < +\infty & \left(c < 1 + \dfrac{1}{m}\right) \end{cases}. \tag{27}$$

This implies that the two parameters $c, m$ cannot be chosen freely with each other to guarantee $C_0 < +\infty$. Specifically, for given $c$, $m$ must be chosen sufficiently small so that

$$c < 1 + \frac{1}{m}. \tag{28}$$

***Remark 2*** Assuming a constant $w$ case ( $w \equiv 1$ as in Yoshioka [29]), then $C_0 < +\infty$ but with the necessary restriction $c = m + 1 > 1$ because $m > 0$. Moreover, (28) combined with $c = m + 1$ requires $m + 1 < 1 + \dfrac{1}{m}$, and hence $c = m + 1 \in (1, 2)$. The real data suggest that both $c < 1$ and $c > 1$ ( as well as $c = 1$ ) are possible [24], and that the constant $w$ case is valid only if $c \in (1, 2)$, implying that large values of $c$ are not covered. This is the primary theoretical justification for considering the time-dependent $w$.

Meaning of the specific coefficient $w$ in (25) is discussed to obtain its insights into fish migration. As shown in (25), behavior of this $w$ is different depending on the sign of $c - (m + 1)$; it is decreasing and vanishes at time $T$ if $c > m + 1$, constant if $c = m + 1$, and increasing and blows up at time $T$ if $c < m + 1$. The constant case assumes that the cost of migration is time-homogeneous and based on **Remark 2**, it is the least interesting case because of the severe limitation of the applicability range. A larger value of $c$ enforces the fish count to more rapidly become 0 owing to the assumed functional form of $h$; the decreasing case therefore implies that the migration is suppressed in the mean near the terminal time and hence near the sunset in our context. By contrast, the increasing case $c < m + 1$ indicates that migration near the sunset is more cost-effective. Detecting the exact origin of the functional shape of $w$ from fish count data only may be difficult, as it results from multiple factors discussed in **Section 3.1**;

however, it can be used for qualitatively categorizing data in each case study. Finally, by employing (28), ranges of parameters $c, m$ in the decreasing and increasing cases can be summarized as follows:

**(Decreasing case)** $$1+m<c<1+\frac{1}{m} \text{ and } 0<m<1, \tag{29}$$

**(Increasing case)** $$c<1+\min\left\{m,\frac{1}{m}\right\} \text{ and } m>0. \tag{30}$$

The following **Proposition 2** shows that for a particular choice of $w$, we can find the condition under which $u_t^*=h_t=\dfrac{c}{T-t}$. This proposition is important because it provides an explicit solution to **Problem O**. Moreover, it provides hints for addressing general cases.

***Proposition 2 (Convergence to Problem O)***

*Assume (25) and (28). Then,* $\Phi_\eta(t,x)\underset{\eta\to+0}{\to}\Phi(t,x)$ *pointwise, and specifically,*

$$\Phi_{\eta_2}(t,x)\le\Phi_{\eta_1}(t,x)\le\Phi(t,x)=A_t x+C_t,\quad 0<t<T,\quad x\ge 0 \tag{31}$$

*with* $A$ *and* $C$ *given in (19) and (20), respectively. Here,* $\eta_1,\eta_2>0$ *are any real numbers such that* $\eta_1\le\eta_2$. *Moreover,* $h_t=\dfrac{c}{T-t}$ *(* $0<t<T$ *) is the optimal control* $u^*$ *of* ***Problem O****, and the corresponding controlled process is given by (1) with* $X_0=X_T=0$.

We also present a generalized version of **Proposition 2** where $w$ is designed using formula (24) such that $h$ satisfies the condition (3).

***Proposition 3 (Generalization of Proposition 2)***

*Assume (24) and (3). Then,* $\Phi_\eta(t,x)\underset{\eta\to+0}{\to}\Phi(t,x)$ *pointwise, and specifically the inequality (31) holds true. Moreover,* $h_t$ *(*$0<t<T$*) is the optimal control* $u^*$ *of* ***Problem O****, and the corresponding controlled process is given by (1) with* $X_0=X_T=0$.

***Remark 3*** The classification of the decreasing and increasing cases (29) and (30) carries over to the generic $h$ that satisfies (3). Indeed, from (24), for $t$ close to $T$, we obtain

$$w_t=O\left((T-t)^{-m}\right)\times e^{-\frac{mc}{m+1}\ln\left(\frac{1}{T-t}\right)}=O\left((T-t)^{\frac{mc}{m+1}-m}\right). \tag{32}$$

***Remark 4*** By inspecting the proofs, **Propositions 1-3** above and **Proposition 4** below carry over to cases where parameter $a,\sigma$ are smooth, bounded, and positive functions of $t\in[0,T]$. This case is discussed in **Section 4**. We conjecture that a similar approach carries over to more complex control problems having singular terminal conditions by exploiting the recent penalization method based on backward SDEs [e.g., 68], which is currently under investigation.

We conclude **Section 3** with the estimates of **Problem O** and **Problem P**. The proof of **Proposition 4** is omitted because it follows from direct calculations.

***Proposition 4***

*Under the assumption of* ***Proposition 3****, for a sufficiently small* $\eta > 0$*, it follows that*

$$\left|A_t - A_{\eta,t}\right| = O(\eta) \quad \text{and} \quad \left|C_t - C_{\eta,t}\right| = O(\eta), \quad 0 < t < T, \tag{33}$$

$$\left|u_t^* - u_{\eta,t}^*\right| = O\left(\eta^{\frac{1}{m}}\right), \quad 0 < t < T, \tag{34}$$

*and*

$$\left|\Phi(t,x) - \Phi_\eta(t,x)\right| = O(\eta) \quad \text{on each compact set in} \quad (0,T)\times(0,+\infty). \tag{35}$$

## 4. Application

### 4.1 Target fish

The target fish in the application of this study is the amphidromous fish species *Ayu*, which has a unique one-year life cycle [69-71]. Each generation starts in autumn, with larvae hatching and migrating downstream to coastal areas, where they overwinter and grow into juveniles. From spring to summer, the juveniles migrate upstream into rivers [72-74], which completes their life cycle in freshwater and uses a lake or dam reservoir instead of the sea.

*Ayu* is a commercially important inland fish species in Japan that constitute almost half of the total production of inland fisheries in the country[1]. The success of each stage (downstream larval migration, winter growth, upstream juvenile migration, and summer growth) of the life history is important in biological, ecological, and fisheries standpoints. Juvenile *Ayu* is known to actively swim and jump during the daytime [23,75], to which the diffusion bridge applies [24]. According to Tsukamoto et al. [76], the origin of the juvenile upstream migration of *Ayu* was hypothesized to be owing to three cascading steps: age and size, endocrinological condition, and psychological processes; however, their detailed mathematical descriptions have not yet been found. Rather than investigating detailed mechanistic models, we employ a diffusion bridge to analyze fish migration from a phenomenological viewpoint. We focus on upstream juvenile migration as a case study owing to data availability, as explained below.

### 4.2 Study site

[1]Situation surrounding inland fisheries and aquaculture, Ministry of Agriculture, Forestry and Fisheries. https://www.jfa.maff.go.jp/j/enoki/attach/pdf/naisuimeninfo-44.pdf (In Japanese. Last accessed on July 28, 2025)

The study site is the Nagara River Estuary Barrage (**Figure 2**), in the Nagara River, Kiso River system, Tokai Region, Japan, where 10-min fish count data of juvenile *Ayu* from late February to late June have been available since 2023. The Nagara River is a major habitat for *Ayu* after the juvenile upstream migration to the middle reaches in Gifu Prefecture [77]. Moreover, the Nagara River and *Ayu* migrating this river have been registered as a "Globally Important Agricultural Heritage System" [78] owing to their strong connection to the circulation of regional human activities, water environment, and fisheries.

An AI-assisted video camera system has been installed at a fishway of the Nagara River Estuary Barrage, where juvenile *Ayu* migrants are counted automatically during the day by the Ibi River, and Nagara River General Management Office[2]. Usually, juvenile migrants of *Ayu* are counted manually (e.g., counted visually or using traps) each day, and sub-daily count data are rare in Japan (see Yoshioka [79] and references therein) possibly because of difficulties in installing advanced observation equipment in natural river environments and their maintenance costs. The Nagara River Estuary Barrage is one of the few observation sites in Japan where sub-daily (more specifically, sub-hourly) fish count data of juvenile *Ayu* have been observed, such that the data can be used for research.

The 10-min fish migration data at the study site are available for the years 2023 (February 22 to June 30) and 2024 (February 26 to June 30) and have been analyzed statistically by Yoshioka [24] (see, Section 4 and Appendix of this literature). This previous study found that the upstream fish migration of *Ayu* at the study site occurred during the daytime and exhibited high fluctuation and intermittency, suggesting the relevance of applying a diffusion bridge. However, these data and model were not studied from the viewpoint of the minimization principle presented in this paper.

[2] NEC Corporation. https://jpn.nec.com/press/201911/20191128_02.html (In Japanese. Last accessed on August 29, 2025)

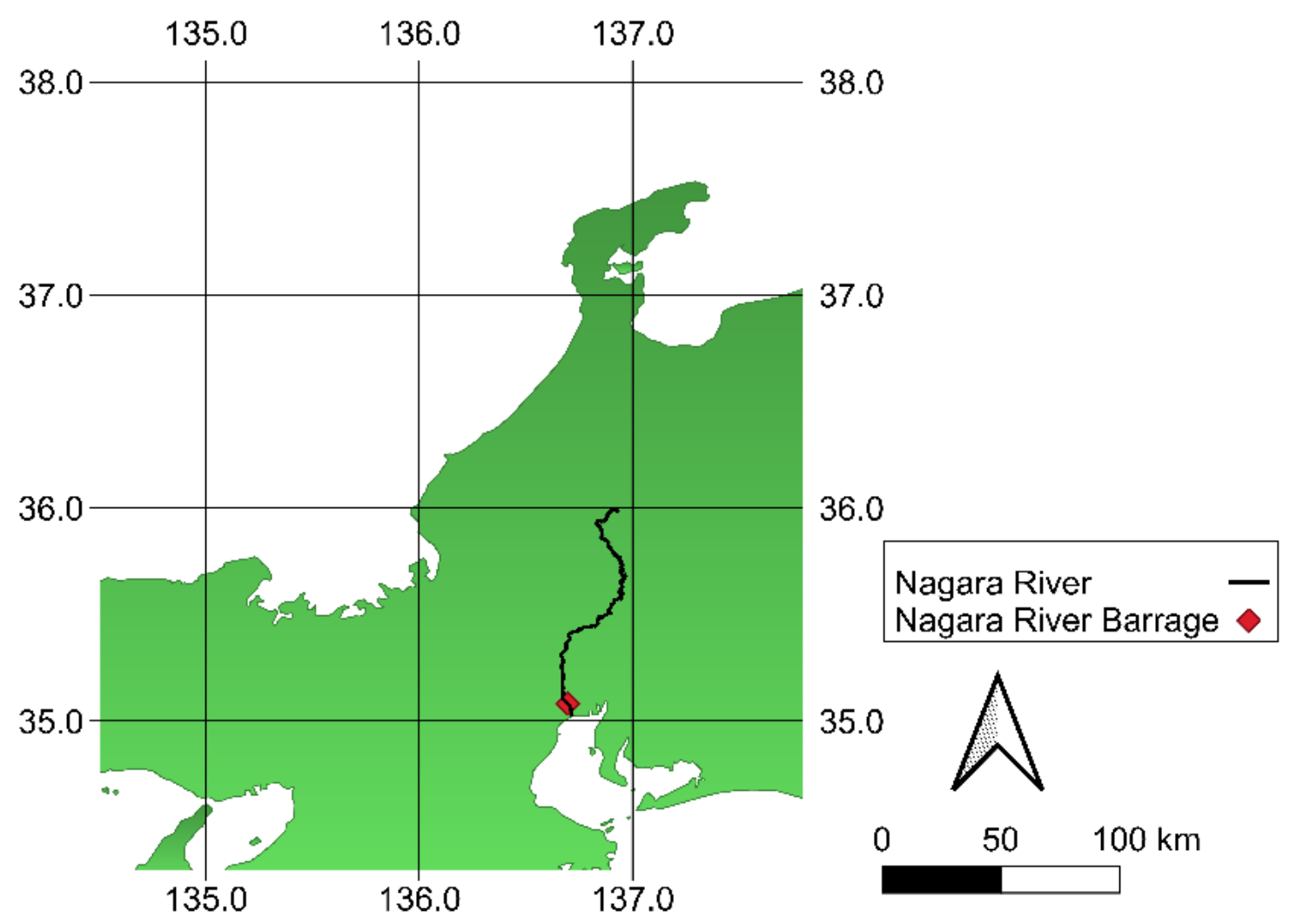

**Figure 2.** Study site.

#### 4.3 Model fitting

##### 4.3.1 Data and correlation

We follow the procedure of Yoshioka [24] with an extended model, as explained below. Each day is normalized to a unit time interval $(0,1)$ and a unit time fish count $X$ to be non-dimensional, where time 0 is sunrise and 1 is sunset. For each day with a positive total daily fish count, the time is normalized by daytime length, the unit-time fish count by the measurement unit (10 min), and the total daily fish count. With this normalization, we regard each path of the daytime fish count as a sample path of the normalized $X$. This data treatment can be at least partly justified if empirical sample paths of fish count are considered uncorrelated between successive days, i.e., between nearest samples. The correlation of the fish count data between successive days is calculated for 2023 and 2024, where each correlation value is considered as a random variable. Taking these correlations would be allowed if daytime lengths between each successive days are close to each other, which is justified because the difference between successive days is less than 10 minutes. **Figure 3** shows the computed correlations between the fish counts on successive days. For both 2023 and 2024, the average correlation is between 0.07 and 0.08 with the corresponding standard deviation approximately 0.17 (**Table 1**). The size of the correlation of the fish count data between successive days is thus considered small.

**Table 1.** Average, standard deviation, and skewness of empirical correlations between successive days in 2023 and 2024.

| | 2023 | 2024 |
|---|---|---|
| Average | 7.43.E-02 | 7.06.E-02 |
| Standard deviation | 1.77.E-01 | 1.66.E-01 |

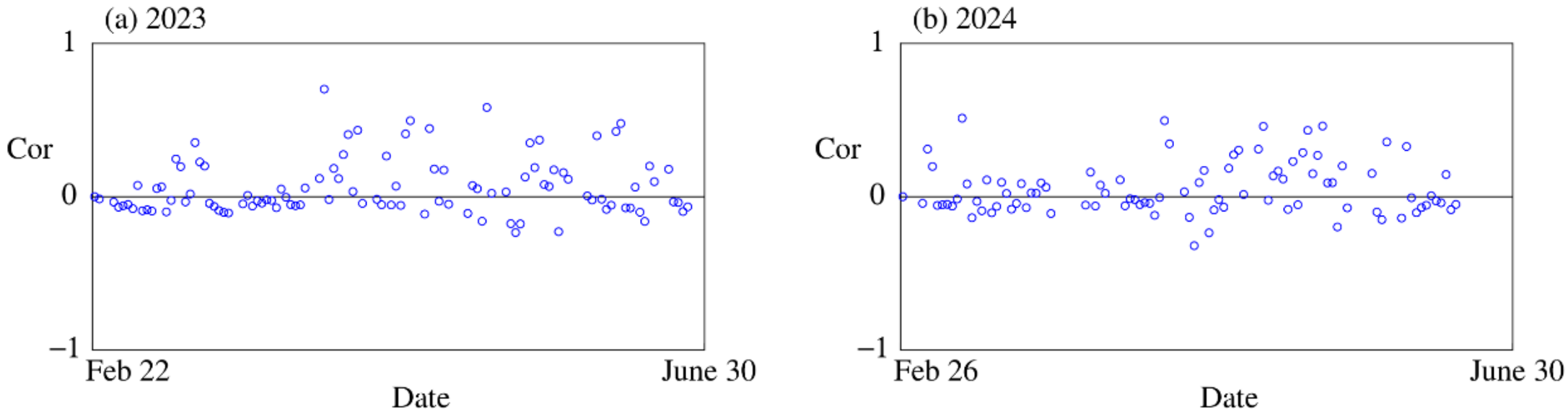

**Figure 3.** Correlation (Cor) of fish counts between each pair of successive days (a) 2023 and (b) 2024. The plots are presented only if the fish counts are available and positive two days in a row.

### 4.3.2 Parameter setting

We fit the parameters of the diffusion bridge to the data from 2023 and 2024 at the study site. We employ a version of the diffusion bridge (1) where the source $a$ and volatility $\sigma$ are given by the time-dependent functions as follows:

$$a_t = a_0 + a_1(1-t) \quad \text{and} \quad \sigma_t = \sqrt{\kappa_0 + \kappa_1(1-t)}, \quad 0 < t < 1 \tag{36}$$

with constants $a_0, a_1$ such that $a_0 \geq 0$ and $a_1 \geq -a_0$, constants $\kappa_0, \kappa_1$ such that $\kappa_0 \geq 0$ and $\kappa_1 \geq -\kappa_0$, and the reversion rate $h$ given following Yoshioka [24] as

$$h_t = \frac{1}{\varepsilon + t} + \frac{1}{1-t} = \frac{1+\varepsilon}{(\varepsilon+t)(1-t)}, \quad 0 < t < 1 \tag{37}$$

with a parameter $\varepsilon > 0$. Models with constant $a, \sigma$ values have been applied to the data of the Nagara River Estuarine Barrage in the previous study, and the current one extends it with a time-dependent source and volatility. The time dependence of $a$ indicates the rate at which individual fishes are supplied from the downstream reach of the observation point, and it is increasing (resp., decreasing) in daytime if $a_1 < 0$ ($a_1 > 0$), respectively. The time dependence of $\sigma$ implies how stochasticity in the fish count evolves in time, and it is increasing (resp., decreasing) in daytime if $\kappa_1 < 0$ ($\kappa_1 > 0$), respectively. The average and variance of $X$ are found in closed forms, as shown in **Appendix B**. The model parameters are identified by the moment-matching method, where $a_0, a_1$ and $\varepsilon$ are identified by a least-squares fitting between the empirical and theoretical $\left(\mathbb{E}[X_t]\right)_{0<t<1}$, and the remaining parameters $\kappa_0, \kappa_1$ are identified by another least-squares fitting between the empirical and theoretical $\left(\sqrt{\mathbb{V}[X_t]}\right)_{0<t<1}$. This fitting method can fully determine the model by exploiting its analytical tractability.

With the proposed control-theoretical framework, the coefficient $w$ associated with the reversion rate $h$ in (37) is expressed as follows, i.e., (24):

$$w_t = \left(\frac{\varepsilon}{\varepsilon + t}\right)^{\frac{m(m+2)}{m+1}} (1-t)^{-\frac{m^2}{m+1}}, \quad 0 \le t < 1, \tag{38}$$

which diverges as $t \to 1$. Moreover, we can choose $m = 1$ as discussed in **Section 3.3**. In particular, the condition $c < 1 + m^{-1}$ is satisfied because in this case $c = 1$ in (3). This corresponds to the increasing case in **Section 3.3** (i.e., (30)) where migration near the sunset is more cost-effective, which is reflected in the identified models with right-weighted average fish counts as shown in **Figures 4-5** in the next subsection.

As indicated in **Remark 1**, the use of a non-constant $w$ is essential in the present case because the specified $h$ does not arise from constant cases and they do not cover $c = 1$. Hence, the diffusion bridge examined in our application is related to the minimization principle only with a non-constant $w$.

### 4.4 Result and discussion

#### 4.4.1 Parameter values

**Table 2** lists the identified parameter values for the 2023 and 2024 data. The root mean-squares error (RMSE) between the empirical and theoretical statistics (i.e., the average and standard deviation) is presented in **Table 2**. **Figures 4-5** compare the empirical and fitted models for the years 2023 and 2024, respectively.

According to **Table 2**, the signs of $a_1$ differ between 2023 (positive) and 2024 (negative), suggesting that the source increases (resp., decreases) during daytime in 2024 (resp., 2023). By contrast, $\kappa_1$ is positive in both 2023 and 2024, implying that volatility decreases over time during the daytime in both years. Regarding the accuracy of the identified models compared to the constant ones in Yoshioka [24], allowing for the time dependence of the parameters $a, \sigma$ improves the accuracy for both average and standard deviation in 2023, while the accuracy of standard deviation is slightly degraded in 2024 by approximately 1%. The results obtained are considered partly owing to the two-step parameter identification procedure that employs the average and standard deviation separately. Nevertheless, the results suggest that the improvement in model accuracy depends on empirical data and is occasionally small, as is the case for 2024. The collection of more data may resolve this issue in the future.

**Figures 4-5** suggest that incorporating the time dependence of the parameters in $a, \sigma$ is not significant in 2024, because the theoretical average and standard deviation peak in the afternoon, as in the constant-parameter case. In particular, the theoretical standard deviations of the two models are difficult to distinguish, as shown in **Figure 5(b)**. In contrast, the peak of average is predicted to occur earlier in the present model than in the constant-parameter model. These results suggest that incorporating time-

dependent parameters does not always work qualitatively and depends on the data available for model fitting.

About the Feller condition to judge volatility regimes for the identified models, we have

$$a_t < \frac{\sigma_t^2}{2} h_t \quad \text{at all} \quad 0 \le t < 1. \tag{39}$$

The identified models are thus in the high-volatility regime for both 2023 and 2024, similar to the constant-parameter models of Yoshioka [24], implying that the intraday fish migration of *Ayu* in the study site is essentially a fluctuating and intermittent phenomenon. Being in the high-volatility regime implies the existence of certain schooling behavior (i.e., migrating with forming clusters) of migrants at the study site; schooling behavior of *Ayu* has been studied in detail through laboratorial experiments [80,81] and certain algebraic model (Yoshioka [82] and references therein). Although we do not have data to directly quantify the schooling behavior in detail, the results obtained are theoretically interesting because they may be able to potentially connect microscopic swimming behavior to macroscopic migration pattern. The time series data do not directly indicate fish school size, but parameters in the SDE as well as coefficient $w$ are considered to affect the schooling behavior of *Ayu* during upstream migration. The observation of school size would be informative for a mechanistic understanding of the diffusion bridge as a conceptual model. Future surveys of fish counts should include predator observations, as schooling in relatively narrow river environments may reduce individual predation risk from birds and terrestrial animals [83].

**Table 2.** Identified parameter values and RMSE for the data in 2023 and 2024. Here, “Constant-parameter model” refers to “Model 2” in Yoshioka [24].

| | 2023 | 2024 |
|---|---|---|
| $a_0$ | 3.375.E-02 | 7.964.E-02 |
| $a_1$ | 5.255.E-02 | -2.707.E-02 |
| $\kappa_0$ | 1.276.E-01 | 2.619.E-01 |
| $\kappa_1$ | 3.853.E-01 | 1.152.E-01 |
| $\varepsilon$ | 1.842.E-01 | 1.370.E-01 |
| RMSE of Average | 4.007E-03 | 4.814E-03 |
| (Constant-parameter model) | (4.283E-03) | (4.870E-03) |
| RMSE of Standard deviation | 1.839E-02 | 1.773E-02 |
| (Constant-parameter model) | (2.089.E-02) | (1.763E-02) |

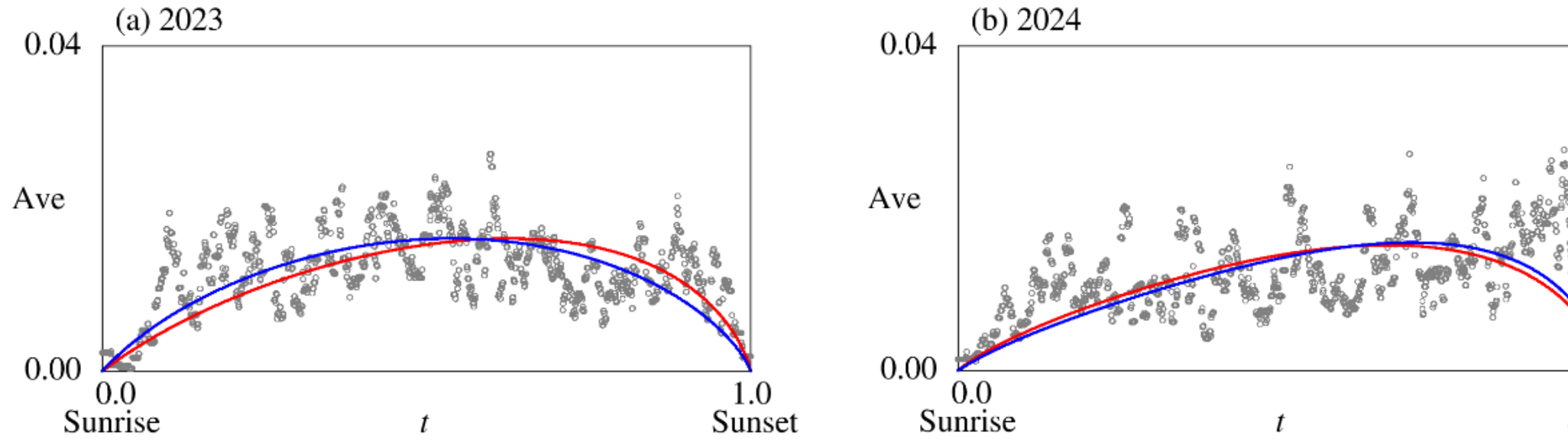

**Figure 4.** Comparison of empirical (circles) and fitted results with the present (blue curve) and constant-parameter (red curve) models in (a) 2023 and (b) 2024.

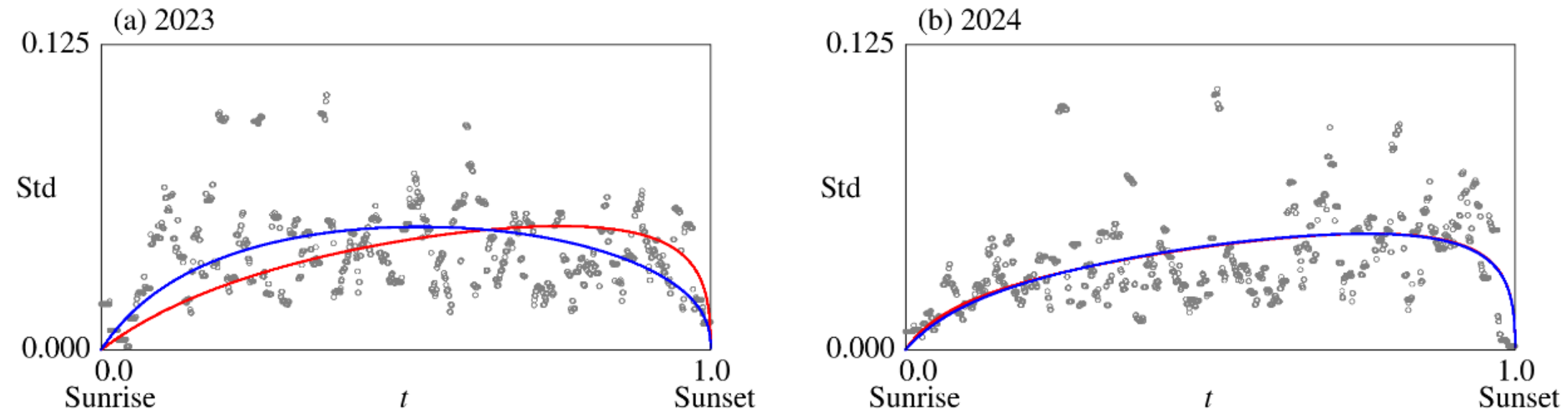

**Figure 5.** Comparison of empirical (circles) and fitted results with the present (blue curve) and constant-parameter (red curve) models in (a) 2023 and (b) 2024.

#### 4.4.2 Cost of migration

We compare the difference between the models in 2023 and 2024 from a viewpoint of the minimization principle. **Figure 6** compares the coefficient $w$ in (38) for the models in 2023 and 2024, where we have chosen $m = 1$. Similarly, **Figure 7** compares $C_t$ for the models in 2023 and 2024. Here, $C_t = \Phi(t,0)$ represents the cost of migration during the time interval $(t,1)$ and is given by

$$C_t = \int_t^1 a_s A_s \mathrm{d}s \, , \quad 0 < t < T \, . \tag{40}$$

**Figure 6** shows that $w_t$ is higher in 2023 than in 2024, implying that the cost of migration is higher in the latter case possibly because the integrand of the objective function $J$ is inversely proportional to $w$. **Figure 7** quantifies the costs of $C$ for the two years, considering both the present and constant-parameter models. For the present model, the cost of migration in 2024 is approximately two to three times of that in 2023. This result is partly owing to the smaller $w$, as explained above, as well as the source of $a$ in 2024 increases in time, whereas that in 2023 decreases. Because the coefficient $A$ assumes large values for small $t$ (i.e., morning) in this case, the unit time cost is further amplified by $a$ during this time period. This finding is consistent with **Figure 4(b)** that the average fish count is higher in the afternoon. In **Figure 7**, a comparison between the present and constant-parameter models implies that the average unit-time fish count of the latter is weighted more (resp., less) in the afternoon and have smaller (resp., larger) unit-time cost values in 2023 (resp., 2024). A monotonic relationship is established between the average unit-time fish count and the unit-time cost of migration.

**Appendix C** discusses sensitivity analysis of the identified models focusing on some manual partial observation schemes for the estimation of daily fish count of *Ayu*. The analysis in that appendix is not directly related to the minimization principle presented in this study, while giving interesting practical results that whole the identified models predict that the intermittent nature of the diffusion bridges implies erroneous estimation results by the conventional manual observation schemes that observe the fish count only a part of daytime. The difference between the identified models, i.e., the difference between their underlying minimization principles, is suggested not to critically affect the results obtained there.

A limitation of the proposed mathematical approach is that the timing of burst events is stochastic, whereas upstream migration in fish populations likely involves decision-making based on current and past environmental cues [84]. However, this may be implicitly represented by the time-dependent coefficients in our SDE, which are linked to the minimization principle. The proposed control formalism essentially depends on time-dependent cues represented by the coefficient $w$ in the objective function; more specifically, this coefficient is possibly given as a function of environmental variables such as water quantity and quality at the study site, which is, unfortunately, currently not available. Nevertheless, our approach could open a new door for deeper understanding of fish migration as a complex biological phenomenon. Collecting more information and data about the intraday upstream migration of *Ayu* along with (temporally fine) environmental data will be an important task towards its better understanding and mathematical formulation.

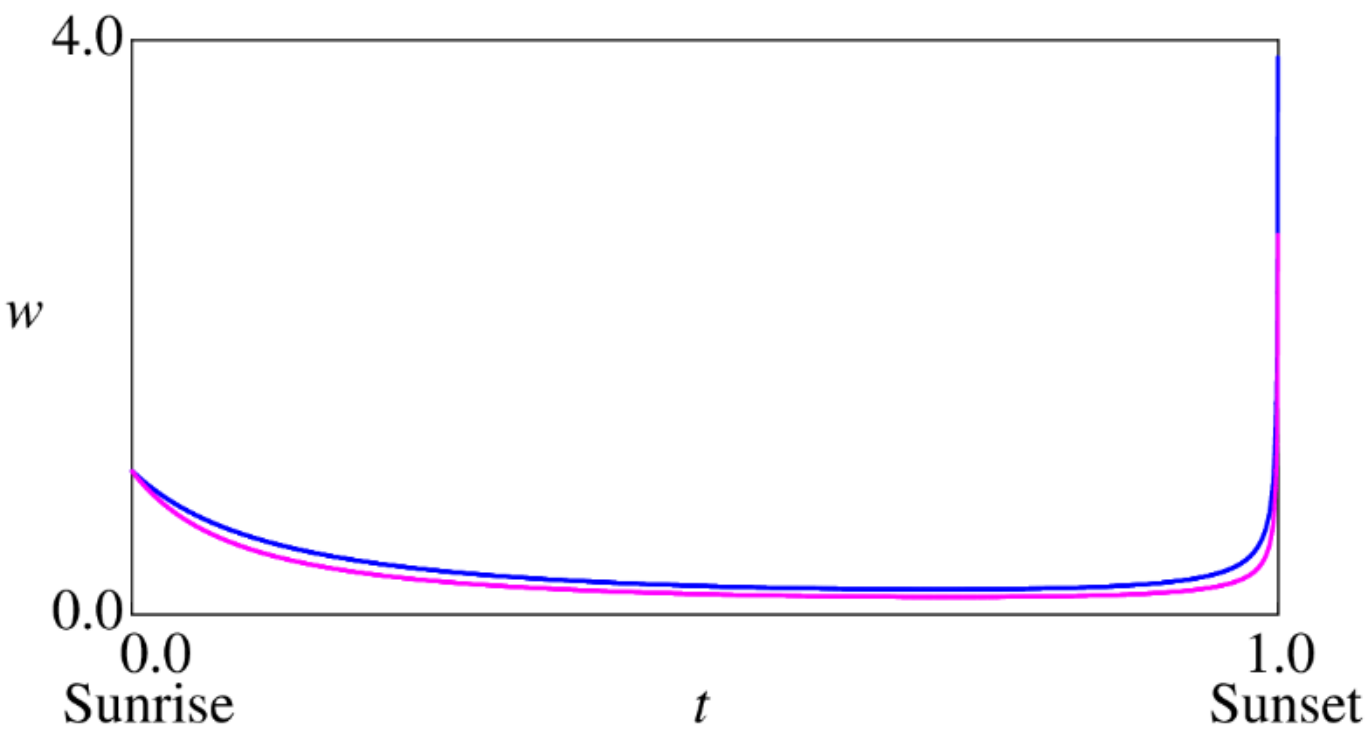

**Figure 6.** Comparison of the coefficient $w$ for 2023 (blue) and 2024 (magenta).

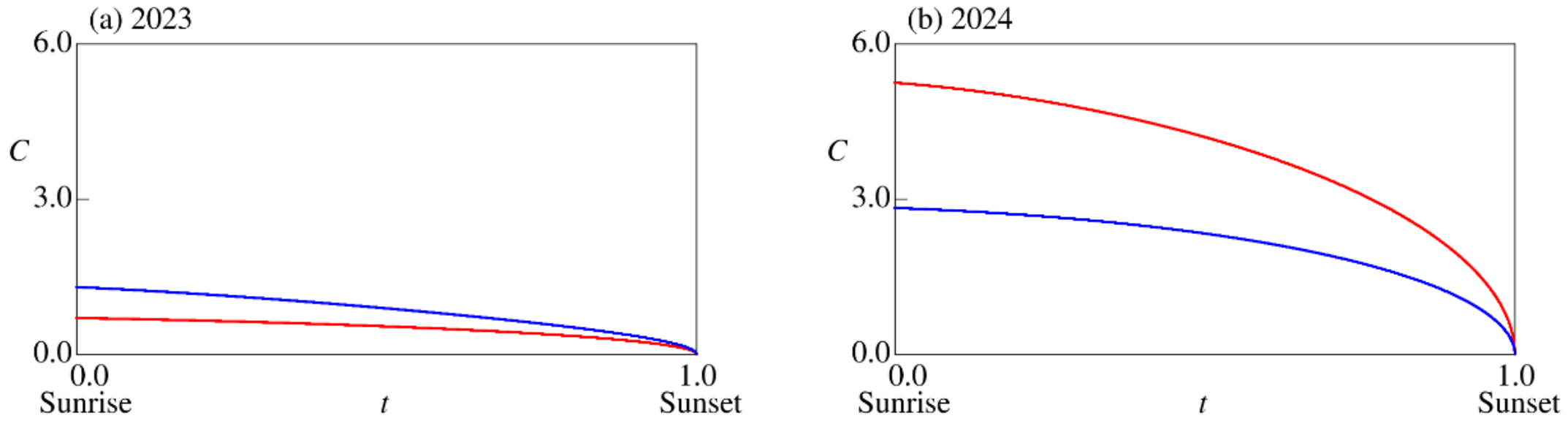

**Figure 7.** Comparison of the cost of migration $C$ for the present model (blue) and constant-parameter model (red) in (a) 2023 and (b) 2024.

## 5. Conclusion

We proposed a minimization principle for intraday fish migration at a fixed observation point. Time-dependent and unbounded coefficients in the diffusion bridge were identified as optimal controls, offering a new theoretical characterization of the fish migration phenomenon. Furthermore, we examined the well-posedness of the original and penalized control problems, establishing the condition under which the penalized problem converges to the original as the penalization vanishes. The application to the juvenile *Ayu* suggested that their diurnal upstream migration can be interpreted through this minimization principle with yearly parameter values.

We applied a diffusion bridge to the diurnal fish migration, although nocturnal migration is also typical for certain fish species, such as eels [85,86]. Finding the proper minimization principle for a diverse form of intraday migration can become feasible with access to high-resolution fish count data (e.g., data in minutes). Categorizing fish species according to their associated energy functionals is essential for advancing our understanding of fish biology and ecology. For greater model flexibility, randomizing initial and terminal conditions may be considered. Thus, sophisticated tools may be used in flow matching [e.g., 87]; however, applying such technique for SDEs with a multiplicative noise is a theoretical challenge, particularly due to nontrivial constraints on boundary conditions. Strengthening the mathematical framework will be critical for further progress in modeling fish migration

The proposed minimization principle can be extended to incorporate both daily and intraday timescales. Although more challenging mathematically, such an extension would offer a more realistic description of fish migration. An important direction for future research is the identifiability of the coefficient $w$ in the objective function from fish count data for each migratory fish species. This topic is currently under investigation by the author. In this paper, we focused on the reasoning about the reversion rate because it is the most mathematically non-trivial ingredient in the model; it determines the singularity of the model. Theoretically, one can also consider the other parameters as control variables, but the associated control problem would be completely different from that studied in this paper and will have to be studied separately. Exploring mechanistic background of the proposed control formulation based on a statistical approach [88] will also be an important topic. We expect that the proposed framework for modeling the migrating fish population would potentially become a building block of an advanced socio-ecological framework involving the population persistence analysis [e.g., 89,90] in future, which is currently under consideration. Approaching the minimization principle from a microscopic, individual-based models will also be an important future topic.

## Appendix

### A. Proofs

***Proof of Proposition 1***

We assume the solution of form (15), and substituting its right-hand side to the HJB equation (12) yields

$$\frac{\mathrm{d}A_{\eta,t}}{\mathrm{d}t}x+\frac{\mathrm{d}C_{\eta,t}}{\mathrm{d}t}+aA_{\eta,t}+x\inf_{v\geq 0}\left\{-vA_{\eta,t}+\frac{1}{w_t}\frac{v^{m+1}}{m+1}\right\}=0,\quad 0<t<T,\quad x\geq 0. \tag{41}$$

Assuming $A_{\eta,t}>0$ for $0<t<T$ in (41) yields

$$\inf_{v\geq 0}\left\{-vA_{\eta,t}+\frac{1}{w_t}\frac{v^{m+1}}{m+1}\right\}=-\frac{m}{m+1}w_t^{\frac{1}{m}}\left(A_{\eta,t}\right)^{\frac{1}{m}+1} \tag{42}$$

with the minimizer

$$\arg\min_{v\geq 0}\left\{-vA_{\eta,t}+\frac{1}{w_t}\frac{v^{m+1}}{m+1}\right\}=w_t^{\frac{1}{m}}\left(A_{\eta,t}\right)^{\frac{1}{m}}, \tag{43}$$

which is identical to the right side of (18). Substituting (42) into (41) yields the system of ordinary differential equations:

$$\frac{\mathrm{d}A_{\eta,t}}{\mathrm{d}t}=\frac{m}{m+1}w_t^{\frac{1}{m}}\left(A_{\eta,t}\right)^{\frac{1}{m}+1},\quad 0<t<T \tag{44}$$

and

$$\frac{\mathrm{d}C_{\eta,t}}{\mathrm{d}t}=-aA_{\eta,t},\quad 0<t<T \tag{45}$$

subject to the terminal conditions $A_{\eta,T}=\frac{1}{\eta}$ and $C_{\eta,T}=0$. The first equation (44) is explicitly solved as (16) and the second one as (17). The explicitly found $A_{\eta,t}$ is positive and bounded because of (14).

The candidate for optimal control is denoted by the right-hand side of (18), with which the controlled SDE becomes

$$\mathrm{d}X_s=\left(a-w_s^{\frac{1}{m}}\left(A_{\eta,s}\right)^{\frac{1}{m}}X_s\right)\mathrm{d}s+\sigma\sqrt{w_s^{\frac{1}{m}}\left(A_{\eta,s}\right)^{\frac{1}{m}}\left|X_s\right|}\mathrm{d}B_s,\quad t<s<T. \tag{46}$$

For any $t\in[0,T)$. We show that the SDE (46) with any initial condition, $X_t=x\geq 0$ admits a unique pathwise non-negative solution for $t<s<T$. To show this, we fix any $t\in[0,T)$. The drift coefficient of SDE (46) is affine in $X_s$, and diffusion involves the square-root $\sqrt{X_s}$ (assuming the absolute value inside the square root turns out to be unnecessary, as shown below). Moreover, the coefficient of $w_s^{\frac{1}{m}}\left(A_{\eta,s}\right)^{\frac{1}{m}}$ is positive and bounded in $(t,T)$. Then, the uniqueness of the pathwise solution to the SDE (46) for each $(t,T-\varepsilon]$ with $\varepsilon>0$ being small follows from Proposition 2.13, Chapter 5 in Karatzas and Shreve [91] or Proposition 1.2.2 in Alfonsi [92], its existence from Theorem 4.22, Chapter 5 in Karatzas and Shreve [91] (we can consider time as another state variable, and $w$ restricted to $[0,T-\varepsilon]$ can be extended to a bounded and continuous function on $\mathbb{R}$), and its non-negativity from an argument similar to that in the last paragraph of p. 5-6 in Alfonsi [92] (we can exploit that the SDE (46) is an affine process with a positive, locally bounded and integrable function $w_s^{\frac{1}{m}}\left(A_{\eta,s}\right)^{\frac{1}{m}}$ in $(t,T)$). For the existence and uniqueness parts, one may also refer to Theorem 6.5 in Chapter 2 of Ma [93].

The assumption (14) suggests that the coefficient $w_s^{\frac{1}{m}}$ is bounded on $[0,T]$ or diverges at most at the rate of $(T-s)^{-\alpha}$ with some $\alpha\in(0,1)$ at the terminal time $T$. In the former case, the unique solution to the SDE (46) can be extended to the full-interval $[0,T]$. In the latter case, we must show that the solution is bounded with probability 1 at $T$. To show this, for $t<s<T$, we have

$$\frac{\mathrm{d}}{\mathrm{d}s}\mathbb{E}[X_s]=a-w_s^{\frac{1}{m}}\left(A_{\eta,s}\right)^{\frac{1}{m}}\mathbb{E}[X_s] \tag{47}$$

and hence

$$\mathbb{E}[X_s]=x\exp\left(-\int_t^s w_u^{\frac{1}{m}}\left(A_{\eta,u}\right)^{\frac{1}{m}}\mathrm{d}u\right)+a\int_t^s\exp\left(-\int_q^s w_u^{\frac{1}{m}}\left(A_{\eta,u}\right)^{\frac{1}{m}}\mathrm{d}u\right)\mathrm{d}q\le x+aT<+\infty . \tag{48}$$

We also have the variance

$$\begin{aligned}\mathbb{V}[X_s]&=\mathbb{E}\left[\left(\mathbb{E}[X_s]-X_s\right)^2\right]\\&=\sigma^2\int_t^s\mathbb{E}\left[X_q\right]w_q^{\frac{1}{m}}\left(A_{\eta,q}\right)^{\frac{1}{m}}\exp\left(-2\int_q^s w_u^{\frac{1}{m}}\left(A_{\eta,u}\right)^{\frac{1}{m}}\mathrm{d}u\right)\mathrm{d}q\\&\le\sigma^2\left(x+aT\right)\int_0^T w_q^{\frac{1}{m}}\left(A_{\eta,q}\right)^{\frac{1}{m}}\mathrm{d}q\\&<+\infty\end{aligned}, \tag{49}$$

with the help of (14) and (16). Moreover, we can rewrite (46) as follows:

$$X_s=\mathbb{E}[X_s]+\sigma\int_t^s\sqrt{w_q^{\frac{1}{m}}\left(A_{\eta,q}\right)^{\frac{1}{m}}X_q}\exp\left(-\int_q^s w_u^{\frac{1}{m}}\left(A_{\eta,u}\right)^{\frac{1}{m}}\mathrm{d}u\right)\mathrm{d}B_q\equiv\mathbb{E}[X_s]+M_s,\quad t<s<T . \tag{50}$$

In view of (48)-(50), $M$ is a square-integrable martingale in each $(t,T-\varepsilon]$ with $\varepsilon>0$ being small. Subsequently, we obtain $\lim_{s\nearrow T}\mathbb{E}[M_s]=0$ and $\lim_{s\nearrow T}\mathbb{E}\left[\left(M_s\right)^2\right]=\lim_{s\nearrow T}\mathbb{V}[X_s]<+\infty$. Therefore, we can extend $(X_s)_{t\le s<T}$ to time and $T$ moreover $X_T<+\infty$ with probability 1. Consequently, the unique solution to the SDE (46) can be extended to the full interval $[0,T]$ under (14).

Finally, the admissibility of the control $w_s^{\frac{1}{m}}\left(A_{\eta,s}\right)^{\frac{1}{m}}$ is verified. To show this, for $0<t<T$ and $x\ge 0$, from (46) we have

$$\begin{aligned}\mathbb{E}\left[\int_t^T \frac{1}{m+1}\frac{1}{w_s}X_s\left(w_s^{\frac{1}{m}}\left(A_{\eta,s}\right)^{\frac{1}{m}}\right)^{m+1}\mathrm{d}s\,\middle|\,X_t=x\right] &= \frac{1}{m+1}\mathbb{E}\left[\int_t^T w_s^{\frac{1}{m}}\left(A_{\eta,s}\right)^{\frac{1}{m}+1}X_s\mathrm{d}s\,\middle|\,X_t=x\right] \\ &= \frac{1}{m+1}\int_t^T w_s^{\frac{1}{m}}\left(A_{\eta,s}\right)^{\frac{1}{m}+1}\mathbb{E}\left[X_s\right]\mathrm{d}s \\ &\leq \frac{x+aT}{m+1}\int_0^T w_s^{\frac{1}{m}}\left(A_{\eta,s}\right)^{\frac{1}{m}+1}\mathrm{d}s \\ &< +\infty\end{aligned} \quad , \tag{51}$$

where we used (14) and the boundedness of $A_\eta$ to obtain the last line. Furthermore, we have

$$\mathbb{E}\left[\frac{1}{\eta}X_T\right] \leq \frac{1}{\eta}\left(x+aT\right) < +\infty \,. \tag{52}$$

Therefore, $w_s^{\frac{1}{m}}\left(A_{\eta,s}\right)^{\frac{1}{m}}$ ( $0<s<T$ ) is an admissible control.

We thus found an explicit classical solution to the HJB equation (12) and proved the unique existence of a pathwise solution to the associated SDE (46). Subsequently, we apply Theorem 3.5.2 of Pham [28] to our penalized control problem because the explicit solution to the HJB equation (12) grows linearly for $x \geq 0$. Here, we note that Theorem 3.5.2 in Pham [28] originally assumed the Lipschitz continuity of diffusion with respect to state variables to guarantee the unique existence of pathwise solutions to the controlled SDE. Although our controlled SDE (46) has a non-Lipschitz diffusion coefficient, this is not problematic because we established the unique existence of its pathwise solution via a different route, as discussed above.

□

***Proof of Proposition 2***

Expectation conditioned on $X_t^{(\cdot)}=\left(X_{\cdot,t}^{(\cdot)}=\right)x\geq 0$ is denoted as $\mathbb{E}^{t,x}$. First, under the assumptions of the proposition, **Proposition 1** holds true. Moreover, we have

$$A_t x + C_t \underset{\eta\to+0}{\to} A_{\eta,t}x + C_{\eta,t} \quad \text{at all} \quad 0\leq t<T \quad \text{and} \quad x\geq 0\,. \tag{53}$$

The next step is to demonstrate the convergence of $\Phi_\eta\left(t,x\right)\to\Phi\left(t,x\right)$ at each point and the ordering property (31). We have

$$A_{\eta,t} = \left(\eta^{\frac{1}{m}} + \frac{1}{m+1}\int_t^T w_s^{\frac{1}{m}}\mathrm{d}s\right)^{-m} \leq \left(\frac{1}{m+1}\int_t^T w_s^{\frac{1}{m}}\mathrm{d}s\right)^{-m} = A_t \quad \text{at all} \quad 0\leq t<T\,, \tag{54}$$

and $A_{\eta,t}$ monotonically converges to $A_t$ from below as $\eta\to+0$, and hence $u^*_{\eta,t}$ to $h_t$. Similarly, we have

$$C_{\eta,t}\leq C_t \quad \text{at all} \quad 0\leq t<T \tag{55}$$

and $C_{\eta,t}$ monotonically converge to $C_t$ from below, as $\eta\to+0$. Therefore, $\Phi_\eta$ is bounded from above by $A_t x + C_t$ irrespective of $\eta>0$.

Now, we show

$$\lim_{\eta\to+0}\mathbb{E}^{t,x}\left[\frac{1}{\eta}X_T^{(u_\eta^*)}\right]\to 0\,, \tag{56}$$

stating that influences of the penalization vanish under the limit $\eta\to+0$ as desired; otherwise, the convergence of **Problems P** to **Problem O** becomes meaningless because there remains a finite or infinite gap between them. By employing (48), an elementary calculation along with (25) yields

$$\int_t^T w_s^{\frac{1}{m}}\mathrm{d}s=\int_t^T\left(1-\frac{s}{T}\right)^{\frac{c}{m+1}-1}\mathrm{d}s=\frac{m+1}{c}T\left(1-\frac{t}{T}\right)^{\frac{c}{m+1}}, \tag{57}$$

$$\begin{aligned}\int_t^T w_u^{\frac{1}{m}}\left(A_{\eta,u}\right)^{\frac{1}{m}}\mathrm{d}u&=\int_t^T\left(1-\frac{u}{T}\right)^{\frac{c}{m+1}-1}\left(\eta^{\frac{1}{m}}+\frac{T}{c}\left(1-\frac{u}{T}\right)^{\frac{c}{m+1}}\right)^{-1}\mathrm{d}u\\&=\int_t^T\frac{\mathrm{d}}{\mathrm{d}u}\left\{-(m+1)\ln\left(\eta^{\frac{1}{m}}+\frac{T}{c}\left(1-\frac{u}{T}\right)^{\frac{c}{m+1}}\right)\right\}\mathrm{d}u\\&=-(m+1)\ln\left(\eta^{\frac{1}{m}}\right)+(m+1)\ln\left(\eta^{\frac{1}{m}}+\frac{T}{c}\left(1-\frac{t}{T}\right)^{\frac{c}{m+1}}\right)\\&=\ln\left(1+\eta^{\frac{-1}{m}}\frac{T}{c}\left(1-\frac{t}{T}\right)^{\frac{c}{m+1}}\right)^{m+1}\end{aligned}, \tag{58}$$

and

$$x\exp\left(-\int_t^T w_u^{\frac{1}{m}}\left(A_{\eta,u}\right)^{\frac{1}{m}}\mathrm{d}u\right)=x\exp\left(-\ln\left(1+\eta^{\frac{-1}{m}}\frac{T}{c}\left(1-\frac{t}{T}\right)^{\frac{c}{m+1}}\right)^{m+1}\right)=x\frac{\eta^{1+\frac{1}{m}}}{\left(\eta^{\frac{1}{m}}+\frac{T}{c}\left(1-\frac{t}{T}\right)^{\frac{c}{m+1}}\right)^{m+1}}\,. \tag{59}$$

Then, for $0\le t<T$, we have

$$\lim_{\eta\to+0}\frac{1}{\eta}x\exp\left(-\int_t^T w_u^{\frac{1}{m}}\left(A_{\eta,u}\right)^{\frac{1}{m}}\mathrm{d}u\right)=x\lim_{\eta\to+0}\frac{\eta^{\frac{1}{m}}}{\left(\eta^{\frac{1}{m}}+\frac{T}{c}\left(1-\frac{t}{T}\right)^{\frac{c}{m+1}}\right)^{m+1}}=0\,. \tag{60}$$

By calculations analogous to (58) and (59), we also have

$$a\int_t^T\exp\left(-\int_q^T w_u^{\frac{1}{m}}\left(A_{\eta,u}\right)^{\frac{1}{m}}\mathrm{d}u\right)\mathrm{d}q=a\int_t^T\frac{\eta^{1+\frac{1}{m}}}{\left(\eta^{\frac{1}{m}}+\frac{T}{c}\left(1-\frac{q}{T}\right)^{\frac{c}{m+1}}\right)^{m+1}}\mathrm{d}q\,. \tag{61}$$

Then, for $0\le t<T$, we have

$$\lim_{\eta \to +0} \frac{1}{\eta} a \int_t^T \exp\left( -\int_q^T w_u^{\frac{1}{m}} \left( A_{\eta,u} \right)^{\frac{1}{m}} \mathrm{d}u \right) \mathrm{d}q = a \lim_{\eta \to +0} \int_t^T \frac{\eta^{\frac{1}{m}}}{\left( \eta^{\frac{1}{m}} + \frac{T}{c}\left(1 - \frac{q}{T}\right)^{\frac{c}{m+1}} \right)^{m+1}} \mathrm{d}q$$

$$\le a \lim_{\eta \to +0} \int_t^T \frac{\eta^{\frac{1}{m}}}{\left( \eta^{\frac{1}{m}} + \frac{T}{c} \right)^{m+1}} \mathrm{d}q \qquad . \tag{62}$$

$$= 0$$

By (62), (60), and

$$\mathbb{E}^{t,x}\left[ X_T^{\left(u_\eta^*\right)} \right] = x \exp\left( -\int_t^T w_u^{\frac{1}{m}} \left( A_{\eta,u} \right)^{\frac{1}{m}} \mathrm{d}u \right) + a \int_t^T \exp\left( -\int_q^T w_u^{\frac{1}{m}} \left( A_{\eta,u} \right)^{\frac{1}{m}} \mathrm{d}u \right) \mathrm{d}q , \tag{63}$$

we obtain the desired result (56). In addition, elementary calculations show that

$$\mathbb{E}^{t,x}\left[ X_T^{(h)} \right] = \mathbb{E}^{t,x}\left[ \left( X_T^{(h)} - \mathbb{E}\left[ X_T^{(h)} \right] \right)^2 \right] = 0 , \tag{64}$$

and hence the control $u = h$ is admissible in both **Problems O and P**.

For each $0 < t < T$ and $x \ge 0$, we have the following convergence result owing to (56) and (64), and the convergence results of $A_{\eta,t}$ and $C_{\eta,t}$:

$$\begin{aligned} \lim_{\eta \to +0} \Phi_\eta(t,x) &= \lim_{\eta \to +0} J_\eta\left(t,x;u_\eta^*\right) \\ &= \lim_{\eta \to +0} \inf_{u \in U} J_\eta(t,x;u) \\ &= \lim_{\eta \to +0} \left( A_{\eta,t} x + C_{\eta,t} \right) \quad . \\ &= A_t x + C_t \\ &= J(t,x;h) \end{aligned} \tag{65}$$

Essentially, the remaining main task is to show

$$A_t x + C_t = J(t,x;h) = \Phi(t,x) , \tag{66}$$

with which $u^* = h$ and the ordering property (31) are obtained.

For each $\eta_1, \eta_2 > 0$ with $\eta_1 \le \eta_2$, we have

$$J_{\eta_2}(t,x;u) \le J_{\eta_1}(t,x;u) \tag{67}$$

for all $0 < t < T$, $x \ge 0$, and $u \in U$, meaning that

$$A_{\eta_2,t} x + C_{\eta_2,t} = \Phi_{\eta_2}(t,x) \le \Phi_{\eta_1}(t,x) = A_{\eta_1,t} x + C_{\eta_1,t} . \tag{68}$$

Therefore, $\Phi_\eta$ is a monotonically decreasing sequence of $\eta > 0$ bounded from above by $A_t x + C_t$ as $\eta \to 0$. Similarly, for each $\eta > 0$, we have

$$J_\eta(t,x;u) \le J(t,x;u) \tag{69}$$

for all $0 < t < T$, $x \ge 0$, and $u \in U$, meaning that

$$\Phi_\eta(t,x) \le \Phi(t,x) . \tag{70}$$

Moreover, due to the suboptimality of $u = h$ in **Problem O**, we have

$$\Phi(t,x) \leq J(t,x;h) \tag{71}$$

for all $0 < t < T$ and $x \geq 0$. Combining (65) and (71) yields

$$\lim_{\eta \to +0} \Phi_\eta(t,x) = \lim_{\eta \to +0} J_\eta(t,x;u_\eta^*) = J(t,x;h) \leq \Phi(t,x) \leq J(t,x;h). \tag{72}$$

Hence, we obtain (66). This implies that $u^* = h$ and (31) hold true.

□

***Proof of Proposition 3***

Expectation conditioned on $X_t^{(\cdot)} = \left(X_{\cdot,t}^{(\cdot)} = \right) x \geq 0$ is denoted as $\mathbb{E}^{t,x}$. The strategy of the proof is essentially the same as that of **Proposition 2**. First, under the assumptions of proposition, **Proposition 1** holds true. Moreover, we have

$$A_t x + C_t \underset{\eta \to +0}{\to} A_{\eta,t} x + C_{\eta,t} \quad \text{at all} \quad 0 < t < T \quad \text{and} \quad x \geq 0. \tag{73}$$

Furthermore, we obtain $A_{\eta,t} \leq A_t$ (resp., $C_{\eta,t} \leq C_0$), and $A_{\eta,t}$ (resp., $C_{\eta,t}$) monotonically converging to $A_t$ (resp., $C_t$) from below as $\eta \to +0$ at all $0 < t < T$. We need to prove

$$\lim_{\eta \to +0} \mathbb{E}^{t,x}\left[\frac{1}{\eta} X_T^{\left(u_\eta^*\right)}\right] \to 0. \tag{74}$$

We assume that this is proven temporally. Subsequently, as in **Proof of Proposition 2**, we have

$$\mathbb{E}^{t,x}\left[X_T^{(h)}\right] = \mathbb{E}^{t,x}\left[\left(X_T^{(h)} - \mathbb{E}^{t,x}\left[X_T^{(h)}\right]\right)^2\right] = 0, \tag{75}$$

The control $u = h$ is admissible in both **Problems O and P**. For each $0 < t < T$ and $x \geq 0$, we obtain

$$\lim_{\eta \to +0} \Phi_\eta(t,x) = \lim_{\eta \to +0} J_\eta(t,x;u_\eta^*) = \lim_{\eta \to +0} \inf_{u \in U} J_\eta(t,x;u) = J(t,x;h). \tag{76}$$

For each $\eta_1, \eta_2 > 0$ with $\eta_1 \leq \eta_2$, we also have

$$J_{\eta_2}(t,x;u) \leq J_{\eta_1}(t,x;u) \leq J(t,x;u) \tag{77}$$

for all $0 < t < T$, $x \geq 0$, and $u \in U$, and hence

$$\Phi_\eta(t,x) \leq \Phi(t,x) \leq J(t,x;h). \tag{78}$$

Consequently, by (76) and (78), we again obtain the relationship (72), and hence the proposition is proven.

The remaining task is to show (75). To show this, we use the formula (63). We have

$$\int_t^T w_u^{\frac{1}{m}} \left(A_{\eta,u}\right)^{\frac{1}{m}} \mathrm{d}u = \int_t^T w_u^{\frac{1}{m}} \left( \eta^{\frac{1}{m}} + \frac{1}{m+1} \int_u^T w_s^{\frac{1}{m}} \mathrm{d}s \right)^{-1} \mathrm{d}u$$

$$\begin{aligned}
&= \int_t^T \frac{\mathrm{d}}{\mathrm{d}u} \left\{ -\left(m+1\right) \ln \left( \eta^{\frac{1}{m}} + \frac{1}{m+1} \int_u^T w_s^{\frac{1}{m}} \mathrm{d}s \right) \right\} \mathrm{d}u \\
&= -\left(m+1\right) \ln \left( \eta^{\frac{1}{m}} \right) + \left(m+1\right) \ln \left( \eta^{\frac{1}{m}} + \frac{1}{m+1} \int_t^T w_s^{\frac{1}{m}} \mathrm{d}s \right) \\
&= \ln \left( 1 + \eta^{\frac{-1}{m}} \frac{1}{m+1} \int_t^T w_s^{\frac{1}{m}} \mathrm{d}s \right)^{m+1}
\end{aligned}, \tag{79}$$

and

$$x \exp\left( -\int_t^T w_u^{\frac{1}{m}} \left(A_{\eta,u}\right)^{\frac{1}{m}} \mathrm{d}u \right) = x \exp\left( -\ln \left( 1 + \eta^{\frac{-1}{m}} \frac{1}{m+1} \int_t^T w_s^{\frac{1}{m}} \mathrm{d}s \right)^{m+1} \right) = x \frac{\eta^{1+\frac{1}{m}}}{\left( \eta^{\frac{1}{m}} + \frac{1}{m+1} \int_t^T w_s^{\frac{1}{m}} \mathrm{d}s \right)^{m+1}} . \tag{80}$$

Then, for $0 \leq t < T$, we have

$$\lim_{\eta \to +0} \frac{1}{\eta} x \exp\left( -\int_t^T w_u^{\frac{1}{m}} \left(A_{\eta,u}\right)^{\frac{1}{m}} \mathrm{d}u \right) = x \lim_{\eta \to +0} \frac{\eta^{\frac{1}{m}}}{\left( \eta^{\frac{1}{m}} + \frac{1}{m+1} \int_t^T w_s^{\frac{1}{m}} \mathrm{d}s \right)^{m+1}} = 0 . \tag{81}$$

By calculations analogous to (79) and (80), we also have

$$a \int_t^T \exp\left( -\int_q^T w_u^{\frac{1}{m}} \left(A_{\eta,u}\right)^{\frac{1}{m}} \mathrm{d}u \right) \mathrm{d}q = a \int_t^T \frac{\eta^{1+\frac{1}{m}}}{\left( \eta^{\frac{1}{m}} + \frac{1}{m+1} \int_q^T w_s^{\frac{1}{m}} \mathrm{d}s \right)^{m+1}} \mathrm{d}q . \tag{82}$$

Hence, for $0 \leq t < T$, we have

$$\lim_{\eta \to +0} \frac{1}{\eta} a \int_t^T \exp\left( -\int_q^T w_u^{\frac{1}{m}} \left(A_{\eta,u}\right)^{\frac{1}{m}} \mathrm{d}u \right) \mathrm{d}q = a \lim_{\eta \to +0} \int_t^T \frac{\eta^{\frac{1}{m}}}{\left( \eta^{\frac{1}{m}} + \frac{1}{m+1} \int_q^T w_s^{\frac{1}{m}} \mathrm{d}s \right)^{m+1}} \mathrm{d}q = 0 \tag{83}$$

due to the classical dominated convergence because the last integrand in (83) is uniformly bounded irrespective of sufficiently small $\eta > 0$. Combining (79)–(83) and (63) yields the desired result (75), and hence completes the proof.

□

**B. Average and variance of the unit-time fish count**

The average and variance of the unit-time fish count $X$ of the diffusion bridge used in **Section 4** are presented below; the average is expressed as follows:

$$\mathbb{E}\left[X_t\right] = a_0 \frac{1-t}{\varepsilon + t} \left\{ \left(1+\varepsilon\right) \ln \left( \frac{1}{1-t} \right) - t \right\} + a_1 \frac{1-t}{\varepsilon + t} \left( \frac{1}{2} t^2 + \varepsilon t \right) \tag{84}$$

and variance by

$$
\begin{aligned}
\mathbb{V}[X_t] = {} & \kappa_0 a_0 (1+\varepsilon)\left(\frac{1-t}{t+\varepsilon}\right)^2 \frac{1}{1-t}\left\{(2+\varepsilon-t)\ln\left(\frac{1}{1-t}\right)-(2+\varepsilon)t\right\} \\
& +\frac{1}{2}\kappa_0 a_1 (1+\varepsilon)\left(\frac{1-t}{\varepsilon+t}\right)^2 \left(t-2(1+\varepsilon)\ln\left(\frac{1}{1-t}\right)+(1+2\varepsilon)\frac{t}{1-t}\right) \\
& +\kappa_1 a_0 (1+\varepsilon)\left(\frac{1-t}{t+\varepsilon}\right)^2 \left(\frac{1}{2}(1+\varepsilon)\left(\ln\left(\frac{1}{1-t}\right)\right)^2+t-\ln\left(\frac{1}{1-t}\right)\right) \\
& +\frac{1}{2}\kappa_1 a_1 (1+\varepsilon)\left(\frac{1-t}{\varepsilon+t}\right)^2 \left(\frac{1}{2}(1-t)(t+3+4\varepsilon)+(2\varepsilon+1)\ln\left(\frac{1}{1-t}\right)-\frac{1}{2}(3+4\varepsilon)\right)
\end{aligned} . \qquad (85)
$$

The formulae (84) and (85) are fully tractable.

### C. A sensitivity analysis of realistic partial observation schemes

We further compare the identified diffusion bridges in **Section 4** to better understand their performance in an application.

In the main text, we studied the 10-min fish count data acquired using a modern video system, but such data are usually not available, and in many cases individuals of *Ayu* have been counted manually; see the literature and reports [e.g., 94,95,96,97,98,99]. Motivated by this issue, the particular focus of this section is manual partial observation schemes of fish count conducted in past studies and surveys, which estimate the total daily fish count by partially observing the fish migration only a part of daytime. Namely, the estimation of total daily fish count is based on its observation during only a part of the whole daytime (e.g., manually counting fish one to several hours). Such approaches are common in the estimation of daily fish count of *Ayu*, which are less costly compared to the whole-day observation scheme; however, from a practical standpoint, it is not always realistic to manually count migrating fish for a long time such as several hours. It is therefore important to understand how much partial observation schemes are accurate provided that the fish count is highly fluctuating and intermittent as for the identified diffusion bridges. In other words, it would be useful to understand how the stochastic nature of the diffusion bridge affects the daily fish count depending on observation schemes.

We formulate partial observation schemes below. As in **Section 4**, daytime is normalized to the unit interval $(0,1)$. A partial observation scheme here estimates the daily fish count $O_l$, which is the observable:

$$
O_l = \frac{1}{l}\int_{l_1}^{l_2} X_s \mathrm{d}s , \qquad (86)
$$

where $0 < l_1 < l_2 < 1$ and $l = l_2 - l_1$. This observation scheme means that the daily fish count of *Ayu* is estimated by observing fish migration during the time interval $(l_1, l_2) \subset (0,1)$ with the length $l$. Intuitively, specifying a longer $l$ would yield a more accurate estimate of the true fish count given by (e.g., Yoshioka [29])

$$
O = \int_0^1 X_s \mathrm{d}s . \qquad (87)
$$

Nevertheless, statistical dependence of $O_l$ on $l_1, l_2$ is not trivial, which we study numerically in this section. Given $O$ and $O_l$, the relative error $R$ to measure accuracy of the former is defined as

$$R = \frac{O_l}{O} - 1, \tag{88}$$

where we assume $O > 0$, which is satisfied at least numerically in the computational cases discussed below. There is no error if $R = 0$ but is not otherwise.

For simplicity, we set $l_1 = 0.5 - 0.5l$ and $l_2 = 0.5 + 0.5l$ so that the center of the observation window is the time 0.5, which approximately corresponds to noon, the middle of daytime. For each computational case presented below, we numerically simulate sample paths of $10^6$ diffusion bridges using the discretization scheme of Yoshioka [24] with the time increment of $10^{-5}$.

**Tables C1-C4** show the computed average, standard deviation, skewness, kurtosis, minimum, and maximum of the relative error $R$ for the identified model in 2023, that in 2024, the constant-parameter model in 2023, and that in 2024. **Figures C1-C3** summarize the average, standard deviation, and coefficient of variation (CV: standard deviation divided by average) as functions of the duration $l$. Average, standard deviation, and CV are suggested to be concave, convex, and convex functions of $l$ for all the computational cases. Computed standard deviation is comparable among the models, while average is not where the constant-parameter model in 2024 is larger than the other three cases. By contrast, the computed CV is the smallest with the constant-parameter model in 2024 and those of the other three cases are comparable. For all the models, the relative error can be made smaller than 10% if $l = 0.9$ that covers 90% of the whole daytime; nevertheless, CV is around 1 even in this case at which the average and fluctuations are comparable. Computed CVs are larger than 1 and average relative errors are larger than 0.1 when $l \leq 0.8$, and hence errors contained in the existing manual partial observation data are considered not negligible. This point is visually implied in the probability density functions of the relative error shown in **Figure C4**, where we find that significant over- and under-estimations are possible when $l$ is small. For small $l$, the PDF has a peak at $R = -1$, suggesting the risk of significant underestimation of the daily fish count.

Finally, the concavity of the average relative error implies that there exists a unique $l$ value between the average relative error and observation cost, the latter being assumed to be proportional to observation duration (e.g., cost proportional to working hours). The implications obtained in this appendix can also be used for designing observation schemes for counting fish considering their operation costs. From a modeling viewpoint, all the models examined give comparable results, and hence the results obtained here would be robust.

**Table C1.** Average, standard deviation, skewness, kurtosis, minimum, and maximum of the relative error $R$ for the identified model in 2023.

| $l$ | Average | Standard deviation | Skewness | Kurtosis | Maximum | Minimum |
|---|---|---|---|---|---|---|
| 0.2 | 2.10.E-01 | 1.15.E+00 | 1.04.E+00 | 1.61.E-01 | 3.96.E+00 | -9.99.E-01 |
| 0.3 | 2.01.E-01 | 8.95.E-01 | 5.16.E-01 | -8.90.E-01 | 2.32.E+00 | -9.99.E-01 |
| 0.4 | 1.90.E-01 | 7.09.E-01 | 7.75.E-02 | -1.22.E+00 | 1.50.E+00 | -9.98.E-01 |
| 0.5 | 1.76.E-01 | 5.60.E-01 | -3.43.E-01 | -1.09.E+00 | 9.99.E-01 | -9.97.E-01 |
| 0.6 | 1.57.E-01 | 4.32.E-01 | -7.93.E-01 | -4.59.E-01 | 6.66.E-01 | -9.93.E-01 |
| 0.7 | 1.33.E-01 | 3.17.E-01 | -1.34.E+00 | 9.98.E-01 | 4.28.E-01 | -9.85.E-01 |
| 0.8 | 1.03.E-01 | 2.09.E-01 | -2.16.E+00 | 4.59.E+00 | 2.50.E-01 | -9.84.E-01 |
| 0.9 | 6.37.E-02 | 1.02.E-01 | -3.90.E+00 | 1.85.E+01 | 1.11.E-01 | -9.65.E-01 |

**Table C2.** Average, standard deviation, skewness, kurtosis, minimum, and maximum of the relative error $R$ for the identified model in 2024.

| $l$ | Average | Standard deviation | Skewness | Kurtosis | Maximum | Minimum |
|---|---|---|---|---|---|---|
| 0.2 | 2.04.E-01 | 1.14.E+00 | 1.06.E+00 | 1.89.E-01 | 3.96.E+00 | -9.99.E-01 |
| 0.3 | 1.96.E-01 | 8.93.E-01 | 5.26.E-01 | -8.78.E-01 | 2.32.E+00 | -9.99.E-01 |
| 0.4 | 1.86.E-01 | 7.10.E-01 | 8.40.E-02 | -1.22.E+00 | 1.50.E+00 | -9.98.E-01 |
| 0.5 | 1.72.E-01 | 5.63.E-01 | -3.41.E-01 | -1.10.E+00 | 9.99.E-01 | -9.96.E-01 |
| 0.6 | 1.55.E-01 | 4.38.E-01 | -7.98.E-01 | -4.71.E-01 | 6.66.E-01 | -9.95.E-01 |
| 0.7 | 1.32.E-01 | 3.25.E-01 | -1.36.E+00 | 1.00.E+00 | 4.28.E-01 | -9.89.E-01 |
| 0.8 | 1.02.E-01 | 2.18.E-01 | -2.20.E+00 | 4.70.E+00 | 2.50.E-01 | -9.90.E-01 |
| 0.9 | 6.32.E-02 | 1.11.E-01 | -4.05.E+00 | 1.95.E+01 | 1.11.E-01 | -9.78.E-01 |

**Table C3.** Average, standard deviation, skewness, kurtosis, minimum, and maximum of the relative error $R$ for the constant-parameter model in 2023.

| $l$ | Average | Standard deviation | Skewness | Kurtosis | Maximum | Minimum |
|---|---|---|---|---|---|---|
| 0.2 | 2.12.E-01 | 1.14.E+00 | 1.03.E+00 | 1.26.E-01 | 3.93.E+00 | -9.99.E-01 |
| 0.3 | 2.05.E-01 | 8.88.E-01 | 5.05.E-01 | -8.96.E-01 | 2.32.E+00 | -9.99.E-01 |
| 0.4 | 1.94.E-01 | 7.04.E-01 | 6.90.E-02 | -1.21.E+00 | 1.50.E+00 | -9.99.E-01 |
| 0.5 | 1.80.E-01 | 5.55.E-01 | -3.53.E-01 | -1.06.E+00 | 9.98.E-01 | -9.97.E-01 |
| 0.6 | 1.61.E-01 | 4.28.E-01 | -8.13.E-01 | -4.02.E-01 | 6.66.E-01 | -9.96.E-01 |
| 0.7 | 1.38.E-01 | 3.12.E-01 | -1.39.E+00 | 1.18.E+00 | 4.28.E-01 | -9.91.E-01 |
| 0.8 | 1.08.E-01 | 2.04.E-01 | -2.27.E+00 | 5.32.E+00 | 2.50.E-01 | -9.92.E-01 |
| 0.9 | 6.65.E-02 | 9.94.E-02 | -4.39.E+00 | 2.42.E+01 | 1.11.E-01 | -9.76.E-01 |

**Table C4.** Average, standard deviation, skewness, kurtosis, minimum, and maximum of the relative error $R$ for the constant-parameter model in 2024.

| $l$ | Average | Standard deviation | Skewness | Kurtosis | Maximum | Minimum |
|---|---|---|---|---|---|---|
| 0.2 | 2.37.E-01 | 1.16.E+00 | 1.02.E+00 | 8.43.E-02 | 3.96.E+00 | -9.99.E-01 |
| 0.3 | 2.27.E-01 | 8.98.E-01 | 4.86.E-01 | -9.27.E-01 | 2.32.E+00 | -9.99.E-01 |
| 0.4 | 2.14.E-01 | 7.10.E-01 | 4.07.E-02 | -1.23.E+00 | 1.50.E+00 | -9.99.E-01 |
| 0.5 | 1.96.E-01 | 5.60.E-01 | -3.90.E-01 | -1.05.E+00 | 9.99.E-01 | -9.96.E-01 |
| 0.6 | 1.74.E-01 | 4.32.E-01 | -8.58.E-01 | -3.47.E-01 | 6.66.E-01 | -9.96.E-01 |
| 0.7 | 1.46.E-01 | 3.18.E-01 | -1.44.E+00 | 1.27.E+00 | 4.28.E-01 | -9.91.E-01 |
| 0.8 | 1.11.E-01 | 2.11.E-01 | -2.31.E+00 | 5.33.E+00 | 2.50.E-01 | -9.92.E-01 |
| 0.9 | 6.64.E-02 | 1.07.E-01 | -4.25.E+00 | 2.16.E+01 | 1.11.E-01 | -9.81.E-01 |

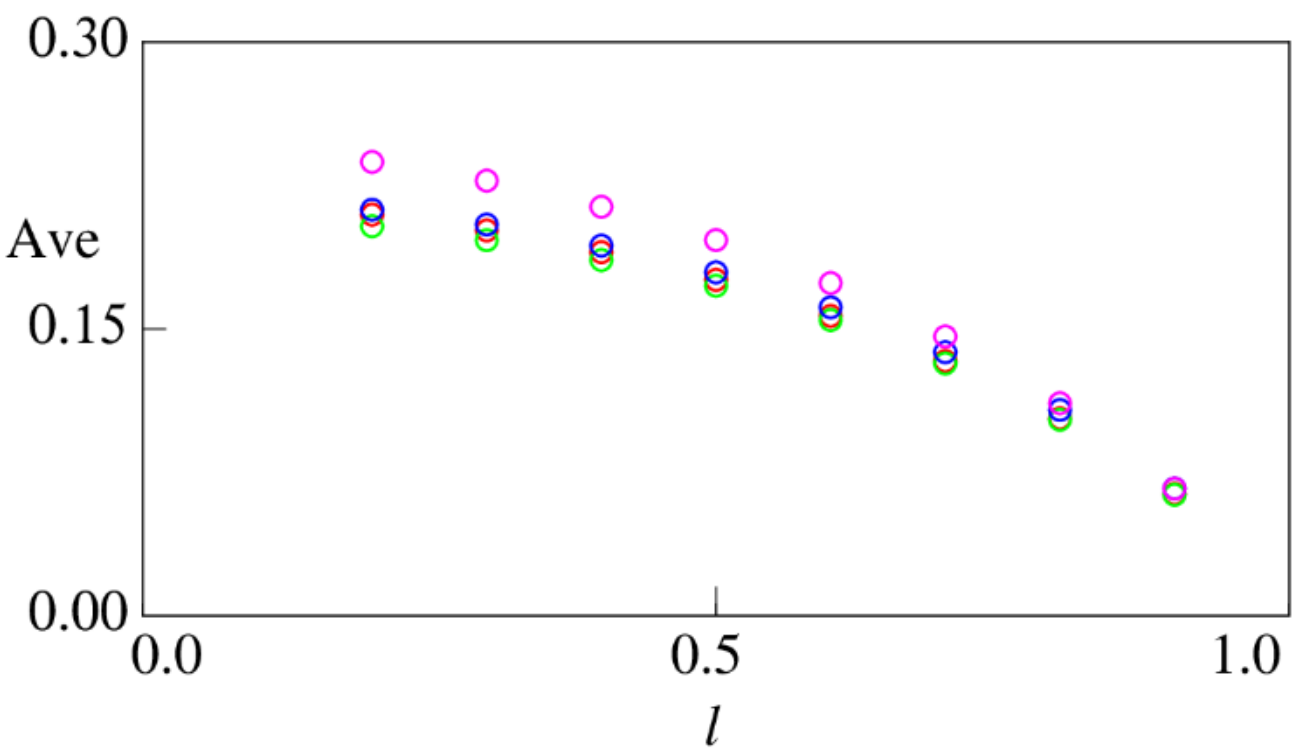

**Figure C1.** Computed average (Ave) of the relative error $R$ as a function of the duration $l$: identified model in 2023 (red), identified model in 2024 (green), constant-parameter model in 2023 (blue), and constant-parameter model in 2024 (magenta).

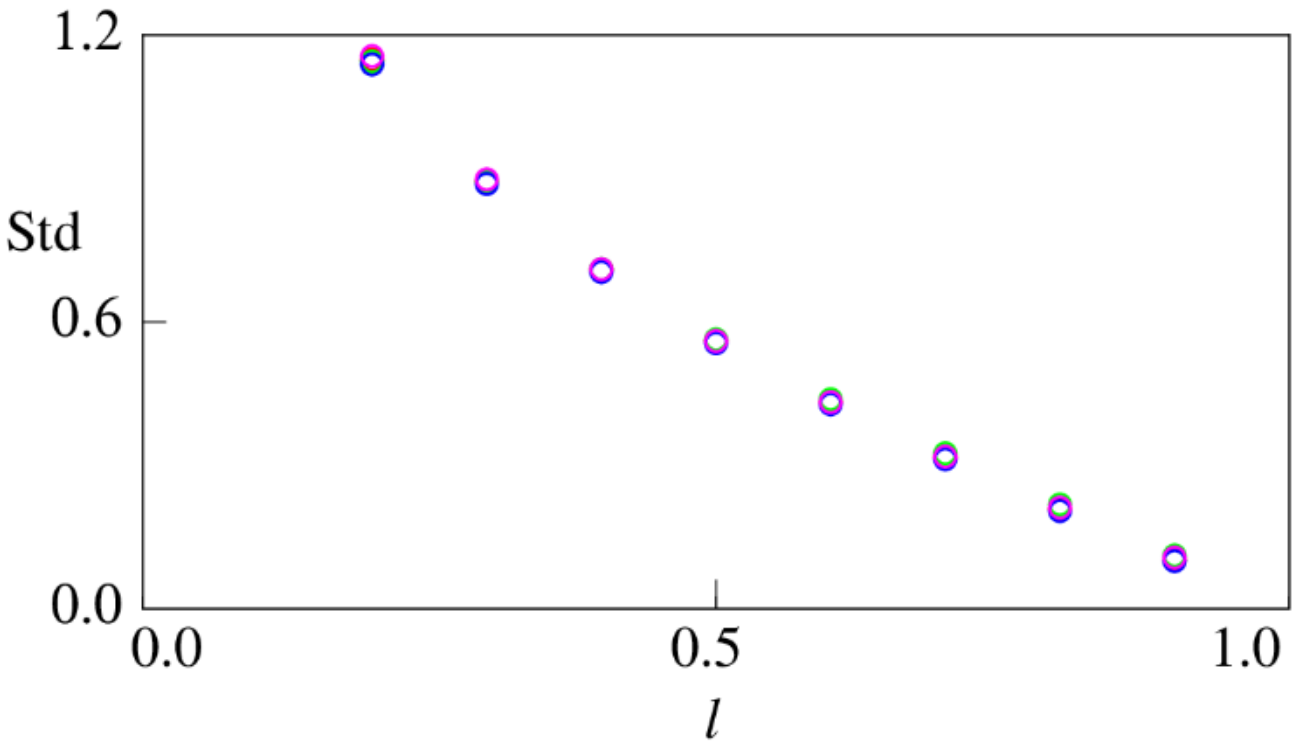

**Figure C2.** Computed standard deviation (Std) of the relative error $R$ as a function of the duration $l$: identified model in 2023 (red), identified model in 2024 (green), constant-parameter model in 2023 (blue), and constant-parameter model in 2024 (magenta).

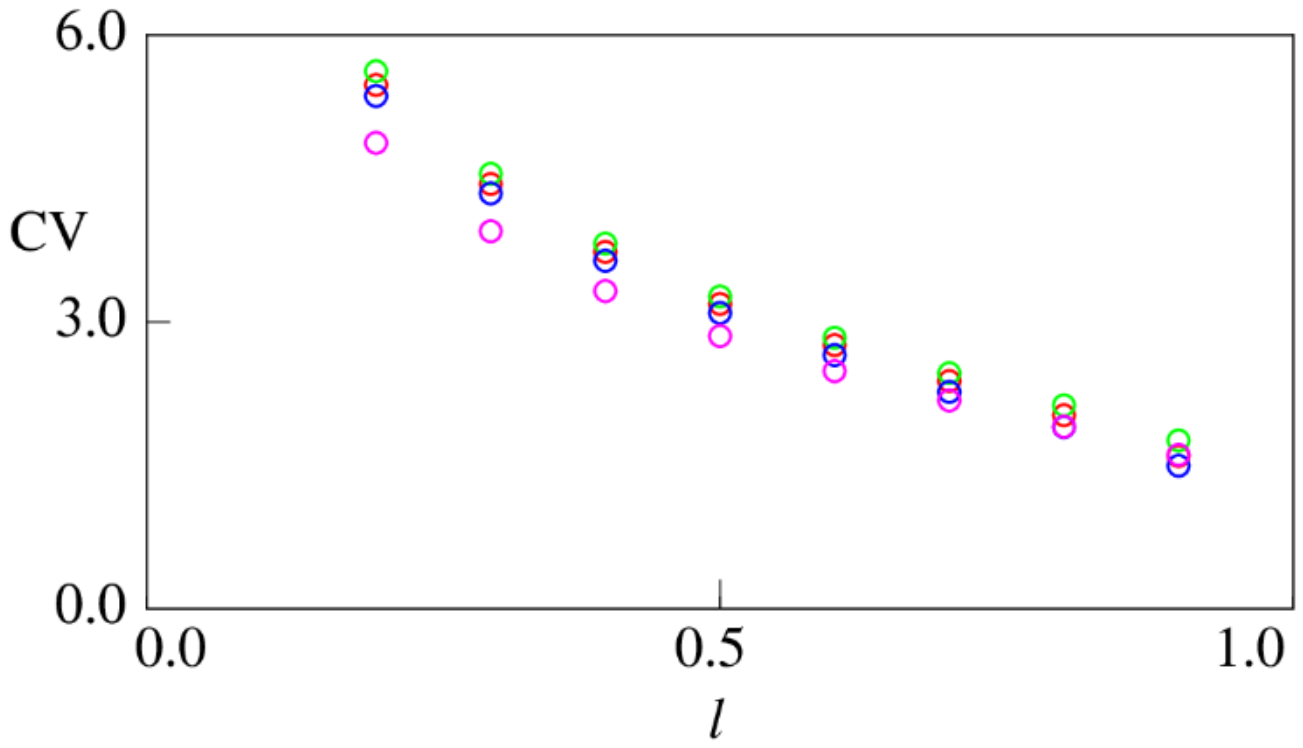

**Figure C3.** Computed coefficient of variation (CV) of the relative error $R$ as a function of the duration $l$: identified model in 2023 (red), identified model in 2024 (green), constant-parameter model in 2023 (blue), and constant-parameter model in 2024 (magenta).

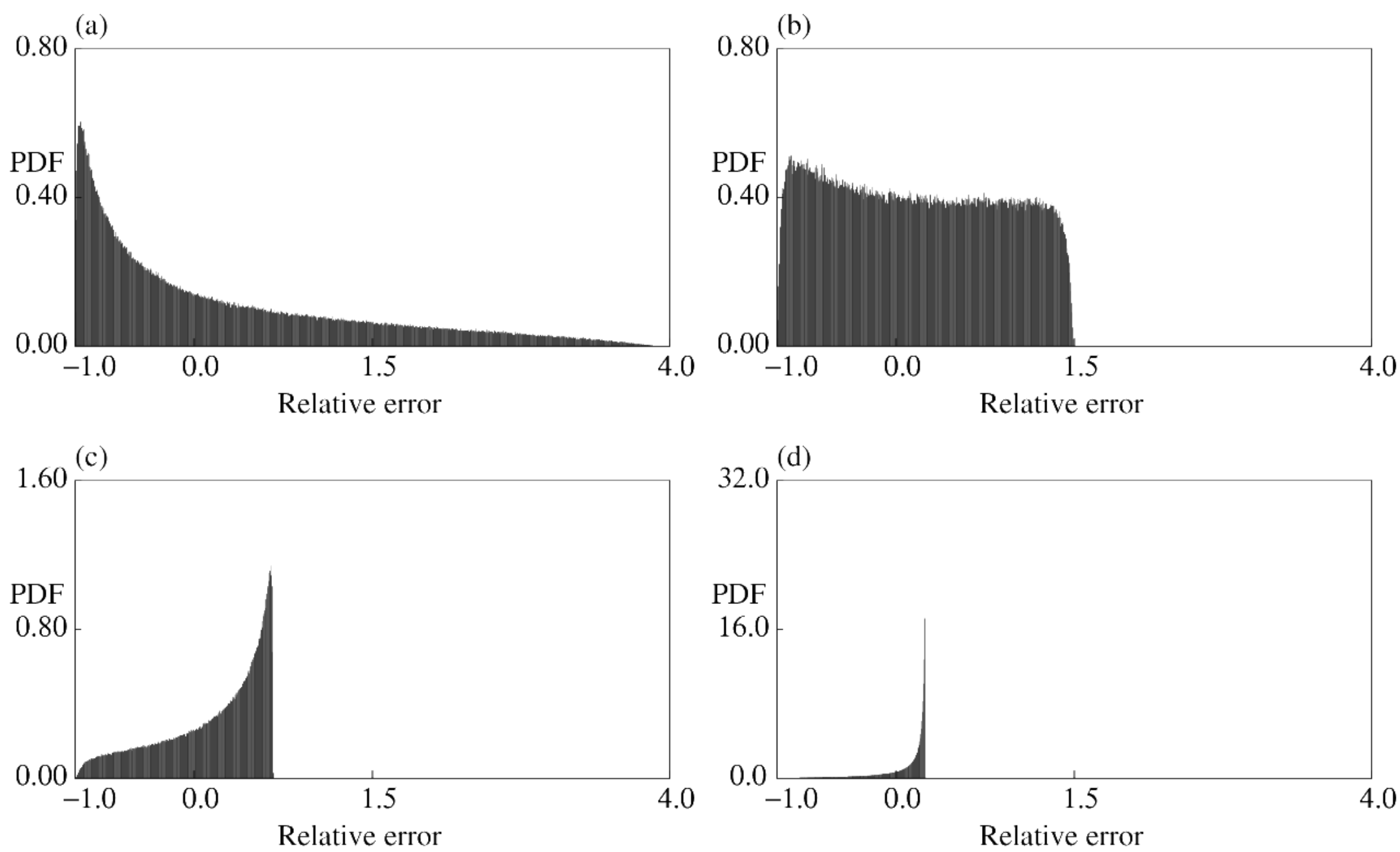

**Figure C4.** Probability density functions (PDFs) of the relative error for the identified model in 2024: (a) $l = 0.2$, (b) $l = 0.4$, (c) $l = 0.6$, and (d) $l = 0.8$.

## References

[1] Alò, D., Lacy, S. N., Castillo, A., Samaniego, H. A., & Marquet, P. A. (2021). The macroecology of fish migration. Global Ecology and Biogeography, 30(1), 99-116. https://doi.org/10.1111/geb.13199

[2] Kurasawa, A., Onishi, Y., Koba, K., Fukushima, K., & Uno, H. (2024). Sequential migrations of diverse fish community provide seasonally prolonged and stable nutrient inputs to a river. Science advances, 10(43), eadq0945. https://doi.org/10.1126/sciadv.adq0945

[3] Hernández-Carrasco, D., Tylianakis, J. M., Lytle, D. A., & Tonkin, J. D. (2025). Ecological and evolutionary consequences of changing seasonality. Science, 388(6750), eads4880. https://doi.org/10.1126/science.ads4880

[4] Kowal, J. L., Haidvogl, G., Funk, A., Schützenhofer, J., Branco, P., Ferreira, M. T., ... & Hein, T. (2025). Over 100 years of longitudinal connectivity changes from the perspective of a migratory fish species. Ecological Indicators, 175, 113436. https://doi.org/10.1016/j.ecolind.2025.113436

[5] Cowx, I. G., Vu, A. V., Hogan, Z., Mallen-Cooper, M., Baumgartner, L. J., Lai, T. Q., ... & Sayer, C. A. (2025). Understanding the threats to fish migration: applying the global swimways concept to the Lower Mekong. Reviews in Fisheries Science & Aquaculture, 33(2), 257-285. https://doi.org/10.1080/23308249.2024.2401018

[6] Booth, E. J., Sandoval-Castillo, J., Attard, C. R., Gilligan, D. M., Unmack, P. J., & Beheregaray, L. B. (2022). Aridification-driven evolution of a migratory fish revealed by niche modelling and coalescence simulations. Journal of Biogeography, 49(9), 1726-1738. https://doi.org/10.1111/jbi.14337

[7] Peluso, L. M., Mateus, L., Penha, J., Bailly, D., Cassemiro, F., Suárez, Y., ... & Lemes, P. (2022). Climate change negative effects on the Neotropical fishery resources may be exacerbated by hydroelectric dams. Science of the Total Environment, 828, 154485. https://doi.org/10.1016/j.scitotenv.2022.154485

[8] Tétard, S., Maire, A., Ovidio, M., Benitez, J. P., Schaeffer, F., Coll, M., & Roy, R. (2025). Multi-year movements of potamodromous cyprinid species within a highly anthropized river assessed using RFID-equipped fishways. Limnologica, 112, 126236. https://doi.org/10.1016/j.limno.2025.126236

[9] Räty, A., Pulkkinen, H., Erkinaro, J., Orell, P., Falkegård, M., & Mäntyniemi, S. (2025). Bayesian species recognition and abundance estimation: unravelling the mysteries of salmonid migration in the Teno River. Canadian Journal of Fisheries and Aquatic Sciences, 82, 1-16. https://doi.org/10.1139/cjfas-2024-0309

[10] Fritts, M., Gibson-Reinemer, D., Appel, D., Lieder, K., Henderson, C., Milde, A., ... & Fritts, A. (2024). Flooding and dam operations facilitate rapid upstream migrations of native and invasive fish species on a regulated large river. Scientific reports, 14(1), 20609. https://doi.org/10.1038/s41598-024-70076-4

[11] Dalton, R. M., Sheppard, J. J., Finn, J. T., Jordaan, A., & Staudinger, M. D. (2022). Phenological variation in spring migration timing of adult Alewife in coastal Massachusetts. Marine and Coastal Fisheries, 14(2), e210198. https://doi.org/10.1002/mcf2.10198

[12] Wansbrough, J., Lokman, P. M., & Closs, G. P. (2024). Local variation in the timing of reproduction and recruitment in a widely distributed diadromous fish. New Zealand journal of marine and freshwater research, 58(2), 238-254. https://doi.org/10.1080/00288330.2023.2212912

[13] Fauchet, L., Laporte, M., Allard, I. C., Moore, J. S., Derôme, N., April, J., & Bernatchez, L. (2025). Monitoring Atlantic Salmon (Salmo salar) smolt migration in a large river system using environmental DNA. Environmental DNA, 7(2), e70087. https://doi.org/10.1002/edn3.70087

[14] Engman, A. C., Kwak, T. J., & Fischer, J. R. (2021). Big runs of little fish: First estimates of run size and exploitation in an amphidromous postlarvae fishery. Canadian Journal of Fisheries and Aquatic Sciences, 78(7), 905-912. https://doi.org/10.1139/cjfas-2020-0093

[15] Griffioen, A. B., van Keeken, O. A., Hamer, A. L., & Winter, H. V. (2022). Passage efficiency and behaviour of sea lampreys (*Petromyzon marinus*, Linnaeus 1758) at a large marine–freshwater barrier. River Research and Applications, 38(5), 906-916. https://doi.org/10.1002/rra.3967

[16] Lilly, J., Honkanen, H. M., McCallum, J. M., Newton, M., Bailey, D. M., & Adams, C. E. (2022). Combining acoustic telemetry with a mechanistic model to investigate characteristics unique to successful Atlantic salmon smolt migrants through a standing body of water. Environmental Biology of Fishes, 105(12), 2045-2063. https://doi.org/10.1007/s10641-021-01172-x

[17] Helminen, J., & Linnansaari, T. (2023). Combining imaging sonar counting and underwater camera species apportioning to estimate the number of Atlantic salmon and striped bass in the Miramichi

River, New Brunswick, Canada. North American Journal of Fisheries Management, 43(3), 743-757. https://doi.org/10.1002/nafm.10889
[18] Hill, C. R., O'Sullivan, A. M., Hogan, J. D., Curry, R. A., Linnansaari, T., & Harrison, P. M. (2025). Poor Upstream Passage for Endangered Atlantic Salmon (*Salmo salar*) at a Trap-and-Transport Fishway in the Wolastoq/Saint John River, New Brunswick. River Research and Applications. Published online. https://doi.org/10.1002/rra.70008
[19] Staples, A., Legett, H. D., Deichmann, J. L., Heggie, K., & Ogburn, M. B. (2023). Automated acoustic detection of river herring (Alewife and Blueback Herring) spawning activity. North American Journal of Fisheries Management, 43(3), 869-881. https://doi.org/10.1002/nafm.10897
[20] Rillahan, C. B., & He, P. (2023). Waiting for the right time and tide: The fine-scale migratory behavior of river herring in two coastal New England streams. Marine and Coastal Fisheries, 15(5), e210273. https://doi.org/10.1002/mcf2.10273
[21] Rato, A. S., Alexandre, C. M., Pedro, S., Mateus, C. S., Pereira, E., Belo, A. F., ... & Almeida, P. R. (2024). New evidence of alternative migration patterns for two Mediterranean potamodromous species. Scientific Reports, 14(1), 23910. https://doi.org/10.1038/s41598-024-74959-4
[22] Rand, P. S., & Fukushima, M. (2014). Estimating the size of the spawning population and evaluating environmental controls on migration for a critically endangered Asian salmonid, Sakhalin taimen. Global Ecology and Conservation, 2, 214-225. https://doi.org/10.1016/j.gecco.2014.09.007
[23] Sato, D., & Seguchi, Y. (2024). Estimation of Ayu-fish migration status and passage time through the Yodo river barrage fishway. Advances in River Engineering, 30, 1-4. In Japanese with English Abstract. https://doi.org/10.11532/river.30.0_1
[24] Yoshioka, H. (2025). CIR bridge for modeling of fish migration on sub-hourly scale. Chaos, Solitons & Fractals, 199(2), 116874. https://doi.org/10.1016/j.chaos.2025.116874
[25] Donahue, M. J., Karnauskas, M., Toews, C., & Paris, C. B. (2015). Location isn't everything: timing of spawning aggregations optimizes larval replenishment. PLoS One, 10(6), e0130694. https://doi.org/10.1371/journal.pone.0130694
[26] Hawley, K. L., Urke, H. A., Kristensen, T., & Haugen, T. O. (2024). Balancing risks and rewards of alternate strategies in the seaward extent, duration and timing of fjord use in contemporary anadromy of brown trout (Salmo trutta). BMC Ecology and Evolution, 24(1), 27. https://doi.org/10.1186/s12862-023-02179-x
[27] Kinnison, M. T., Unwin, M. J., Hendry, A. P., & Quinn, T. P. (2001). Migratory costs and the evolution of egg size and number in introduced and indigenous salmon populations. Evolution, 55(8), 1656-1667. https://doi.org/10.1111/j.0014-3820.2001.tb00685.x
[28] Pham, H. (2009). Continuous-time stochastic control and optimization with financial applications. Springer Berlin, Heidelberg.
[29] Yoshioka, H. (2025-2) Two Issues in Modelling Fish Migration. Preprint. https://www.researchgate.net/publication/395130366_Two_Issues_in_Modelling_Fish_Migration
[30] Frølich, E. F. (2023). Copuling population dynamics and diel migration patterns. Theoretical Population Biology, 151, 19-27. https://doi.org/10.1016/j.tpb.2023.03.004
[31] Mazuryn, M., & Thygesen, U. H. (2023). Mean field games for diel vertical migration with diffusion. Bulletin of Mathematical Biology, 85(6), 49. https://doi.org/10.1007/s11538-023-01154-3
[32] Meyer, A. D., Hastings, A., & Largier, J. L. (2024). Making Your Own Luck: Weak Vertical Swimming Improves Dispersal Success for Coastal Marine Larvae. Bulletin of Mathematical Biology, 86(3), 23. https://doi.org/10.1007/s11538-023-01252-2
[33] Frølich, E. F., & Thygesen, U. H. (2022). Solving multispecies population games in continuous space and time. Theoretical Population Biology, 146, 36-45. https://doi.org/10.1016/j.tpb.2022.06.002
[34] Shaw, A. K., Torstenson, M., Craft, M. E., & Binning, S. A. (2023). Gaps in modelling animal migration with evolutionary game theory: infection can favour the loss of migration. Philosophical Transactions of the Royal Society B, 378(1876), 20210506. https://doi.org/10.1098/rstb.2021.0506
[35] Yoshioka, H., & Yaegashi, Y. (2018). An optimal stopping approach for onset of fish migration. Theory in Biosciences, 137(2), 99-116. https://doi.org/10.1007/s12064-018-0263-8
[36] Yoshioka, H. (2019). A stochastic differential game approach toward animal migration. Theory in Biosciences, 138(2), 277-303. https://doi.org/10.1007/s12064-019-00292-4
[37] Movilla Miangolarra, O., Eldesoukey, A., Movilla Miangolarra, A., & Georgiou, T. T. (2025). Maximum entropy inference of reaction–diffusion models. The Journal of Chemical Physics, 162(19). https://doi.org/10.1063/5.0256659

[38] Oishi, K., Hashizume, Y., Nakao, T., & Kashima, K. (2025). Extraction of implicit field cost via inverse optimal Schrödinger bridge. SICE Journal of Control, Measurement, and System Integration, 18(1), 2490332. https://doi.org/10.1080/18824889.2025.2490332

[39] Park, B., Choi, J., Lim, S., & Lee, J. (2024). Stochastic optimal control for diffusion bridges in function spaces. Advances in Neural Information Processing Systems, 37, 28745-28771. https://proceedings.neurips.cc/paper_files/paper/2024/hash/328b81881da145412f2bc56c998dfb6a-Abstract-Conference.html

[40] Reddy, G. (2025). Dynamic landscapes and statistical limits on growth during cell fate specification. Physical Review E, 111(6), 064418. https://doi.org/10.1103/63d2-4wq6

[41] Zhang, K., Zhu, J., Kong, D., & Zhang, Z. (2024-1). Modeling single cell trajectory using forward-backward stochastic differential equations. PLOS Computational Biology, 20(4), e1012015. https://doi.org/10.1371/journal.pcbi.1012015

[42] Zhang, P., Gao, T., Guo, J., & Duan, J. (2025). Action functional as an early warning indicator in the space of probability measures via Schrödinger bridge. Quantitative Biology, 13(3), e86. https://doi.org/10.1002/qub2.86

[43] Garg, J., Zhang, X., & Zhou, Q. (2024). Soft-constrained Schrödinger Bridge: a stochastic control approach. In International Conference on Artificial Intelligence and Statistics (pp. 4429-4437). PMLR. https://proceedings.mlr.press/v238/garg24a.html

[44] Yoshioka, H., & Yamazaki, K. (2023). A jump Ornstein–Uhlenbeck bridge based on energy-optimal control and its self-exciting extension. IEEE Control Systems Letters, 7, 1536-1541. https://doi.org/10.1109/LCSYS.2023.3271422

[45] Duffie, D., Filipović, D., & Schachermayer, W. (2003). Affine processes and applications in finance. The Annals of Applied Probability, 13(3), 984. https://doi.org/10.1214/aoap/1060202833

[46] Abi Jaber, E. (2024). Simulation of square-root processes made simple: applications to the Heston model. Preprint. https://arxiv.org/abs/2412.11264

[47] Desmettre, S., Leobacher, G., & Rogers, L. C. G. (2021). Change of drift in one-dimensional diffusions. Finance and Stochastics, 25(2), 359-381. https://doi.org/10.1007/s00780-021-00451-w

[48] Mishura, Y., & Yurchenko-Tytarenko, A. (2023). Standard and fractional reflected Ornstein–Uhlenbeck processes as the limits of square roots of Cox–Ingersoll–Ross processes. Stochastics, 95(1), 99-117. https://doi.org/10.1080/17442508.2022.2047188

[49] Mishura, Y., Pilipenko, A., & Yurchenko-Tytarenko, A. (2024). Low-dimensional Cox-Ingersoll-Ross process. Stochastics, 96(5), 1530-1550. https://doi.org/10.1080/17442508.2023.2300291

[50] Fukushima, M., & Rand, P. S. (2023). Individual variation in spawning migration timing in a salmonid fish—Exploring roles of environmental and social cues. Ecology and Evolution, 13(5), e10101. https://doi.org/10.1002/ece3.10101

[51] Knudsen, T. E., MacKenzie, B. R., Thygesen, U. H., & Mariani, P. (2025). Evolution and stability of social learning in animal migration. Movement Ecology, 13(1), 43. https://doi.org/10.1186/s40462-025-00564-3

[52] Massie, J. A., Santos, R. O., Rezek, R. J., James, W. R., Viadero, N. M., Boucek, R. E., ... & Rehage, J. S. (2022). Primed and cued: long-term acoustic telemetry links interannual and seasonal variations in freshwater flows to the spawning migrations of Common Snook in the Florida Everglades. Movement ecology, 10(1), 48. https://doi.org/10.1186/s40462-022-00350-5

[53] Ciraane, U. D., Sonny, D., Lerquet, M., Archambeau, P., Dewals, B., Pirotton, M., ... & Erpicum, S. (2025). Impact of anthropized river hydrodynamic conditions on the downstream migration of Atlantic salmon smolts. Science of The Total Environment, 982, 179556. https://doi.org/10.1016/j.scitotenv.2025.179556

[54] Elings, J., Bruneel, S., Pauwels, I. S., Schneider, M., Kopecki, I., Coeck, J., ... & Goethals, P. L. (2024). Finding navigation cues near fishways. Biological Reviews, 99(1), 313-327. https://doi.org/10.1111/brv.13023

[55] Liu, S., Jian, Y., Li, P., Liang, R., Chen, X., Qin, Y., ... & Li, K. (2023). Optimization schemes to significantly improve the upstream migration of fish: A case study in the lower Yangtze River basin. Ecological Engineering, 186, 106838. https://doi.org/10.1016/j.ecoleng.2022.106838

[56] Lombardo, S. M., Buckel, J. A., Hain, E. F., Griffith, E. H., & White, H. (2020). Evidence for temperature-dependent shifts in spawning times of anadromous alewife (Alosa pseudoharengus) and blueback herring (Alosa aestivalis). Canadian Journal of Fisheries and Aquatic Sciences, 77(4), 741-751. https://doi.org/10.1139/cjfas-2019-014

[57] Wanjari, U. R., Mukherjee, A. G., Gopalakrishnan, A. V., Murali, R., Kannampuzha, S., & Prabakaran, D. S. (2023). Factors Affecting Fish Migration. In: Soni, R., Suyal, D.C., Morales-Oyervides, L., Sungh Chauhan, J. (eds) Current Status of Fresh Water Microbiology. Springer, Singapore. https://doi.org/10.1007/978-981-99-5018-8_20

[58] Crawford, R. M., Gee, E. M., Hicks, B. J., & Franklin, P. A. (2025). Group swimming significantly decreases time to passage success for a galaxiid species. Journal of Fish Biology. Online published. https://doi.org/10.1111/jfb.70040

[59] Griffin, K. R., Holbrook, C. M., Zielinski, D. P., Cahill, C. L., & Wagner, C. M. (2025). Not all who meander are lost: migrating sea lamprey follow river thalwegs to facilitate safe and efficient passage upstream. Journal of Experimental Biology, 228(4), JEB249539. https://doi.org/10.1242/jeb.249539

[60] Zhang, Y., Ko, H., Calicchia, M. A., Ni, R., & Lauder, G. V. (2024-2). Collective movement of schooling fish reduces the costs of locomotion in turbulent conditions. PLoS biology, 22(6), e3002501. https://doi.org/10.1371/journal.pbio.3002501

[61] Martin, B. T., Nisbet, R. M., Pike, A., Michel, C. J., & Danner, E. M. (2015). Sport science for salmon and other species: ecological consequences of metabolic power constraints. Ecology letters, 18(6), 535-544. https://doi.org/10.1111/ele.12433

[62] Hartono, A. D., Nguyen, L. T. H., & Tạ, T. V. (2024). A stochastic differential equation model for predator-avoidance fish schooling. Mathematical Biosciences, 367, 109112. https://doi.org/10.1016/j.mbs.2023.109112

[63] Zhang, Y., & Lauder, G. V. (2025). Physics and physiology of fish collective movement. Newton, 1(1). https://doi.org/10.1016/j.newton.2025.100021

[64] Kumada, N., Arima, T., Tsuboi, J. I., Ashizawa, A., & Fujioka, M. (2013). The multi-scale aggregative response of cormorants to the mass stocking of fish in rivers. Fisheries research, 137, 81-87. https://doi.org/10.1016/j.fishres.2012.09.005

[65] Takai, N., Kawabe, K., Togura, K., Kawasaki, K., & Kuwae, T. (2018). The seasonal trophic link between Great Cormorant Phalacrocorax carbo and ayu *Plecoglossus altivelis altivelis* reared for mass release. Ecological research, 33(5), 935-948. https://doi.org/10.2326/osj.22.183

[66] Säterberg, T., Jacobson, P., Ovegård, M., Rask, J., Östergren, J., Jepsen, N., & Florin, A. B. (2023). Species-and origin-specific susceptibility to bird predation among juvenile salmonids. Ecosphere, 14(12), e4724. https://doi.org/10.1002/ecs2.4724

[67] Sortland, L. K., Wightman, G., Flávio, H., Aarestrup, K., & Roche, W. (2024). A physical bottleneck increases predation on Atlantic salmon smolts during seaward migration in an Irish Index River. Fisheries Management and Ecology, e12779. https://doi.org/10.1111/fme.12779

[68] Cacitti-Holland, D., Denis, L., & Popier, A. (2025). Continuity problem for BSDE and IPDE with singular terminal condition. Journal of Mathematical Analysis and Applications, 543(1), 128845. https://doi.org/10.1016/j.jmaa.2024.128845

[69] Murase, I., & Iguchi, K. I. (2022). High growth performance in the early ontogeny of an amphidromous fish, Ayu *Plecoglossus altivelis altivelis*, promoted survival during a disastrous river spate. Fisheries Management and Ecology, 29(3), 224-232. https://doi.org/10.1111/fme.12524

[70] Yoshida, K., Yajima, H., Islam, M. T., & Pan, S. (2024). Assessment of spawning habitat suitability for Amphidromous Ayu (*Plecoglossus altivelis*) in tidal Asahi River sections in Japan: Implications for conservation and restoration. River Research and Applications, 40(8), 1497-1511. https://doi.org/10.1002/rra.4329C

[71] Watanabe, S., Iida, M., Lord, C., Keith, P., & Tsukamoto, K. (2014). Tropical and temperate freshwater amphidromy: a comparison between life history characteristics of Sicydiinae, ayu, sculpins and galaxiids. Reviews in Fish Biology and Fisheries, 24(1), 1-14. https://doi.org/10.1007/s11160-013-9316-8

[72] Kaewsangk, K., Hayashizaki, K. I., Asahida, T., Nemoto, T., & Ida, H. (2001). Clarification of the origin of landlocked ayu, Plecoglossus altivelis, populations in the Kasumigaura Lake system, Ibaraki Prefecture, Japan. Fisheries science, 67(6), 1175-1177. https://doi.org/10.1046/j.1444-2906.2001.00377.x

[73] Khatun, D., Tanaka, T., Horinouchi, M., Yoshioka, H., & Aranishi, F. (2025). allometric growth and condition factors throughout an annual life history of landlocked ayu *Plecoglossus altivelis altivelis* in Haidzuka Dam Reservoir. Lakes & Reservoirs: Research & Management, 30(1), e70011. https://doi.org/10.1111/lre.70011

[74] Tsuji, S., Shibata, N., Sawada, H., & Watanabe, K. (2023). Differences in the genetic structure between and within two landlocked Ayu groups with different migration patterns in Lake Biwa

revealed by environmental DNA analysis. Environmental DNA, 5(5), 894-905. https://doi.org/10.1002/edn3.345
[75] Uchida, K., Tsukamoto, K., & Kajihara, T. (1990). Effects of environmental factors on jumping behaviour of the juvenile ayu *Plecoglossus altiveli*s with special reference to their upstream migration. Nippon Suisan Gakkaishi, 56(9), 1393-1399. DOI https://doi.org/10.2331/suisan.56.1393
[76] Tsukamoto, K, Miller, M. J., Kotake, A., Aoyama, J., & Uchida, K. (2009). The origin of fish migration: the random escapement hypothesis. In American Fisheries Society Symposium (Vol. 69, pp. 45-61), American Fisheries Society.
[77] Nagai, S., Saitoh, T., & Nagayama, S. (2024). Monitoring the phenology of ayu and satsukimasu trout through wholesale market statistics. Japanese Journal of Biometeorology, 61(1), 19-31. In Japanese with English Abstract. https://doi.org/10.11227/seikisho.61.19
[78] Harada, M., Shiozawa, S., Suzuki, M., & Nagayama, S. (2024). Physical habitat assessment for Ayu spawning habitat focusing on the riverbed environment. Japanese Journal of JSCE, 80(16), 23-16105. https://doi.org/10.2208/jscejj.23-16105
[79] Yoshioka, H. (2025-3). Superposition of interacting stochastic processes with memory and its application to migrating fish counts. Chaos, Solitons & Fractals, 192, 115911. https://doi.org/10.1016/j.chaos.2024.115911
[80] Niizato, T., Sakamoto, K., Mototake, Y. I., Murakami, H., & Tomaru, T. (2024). Information structure of heterogeneous criticality in a fish school. Scientific Reports, 14(1), 29758. https://doi.org/10.1038/s41598-024-79232-2
[81] Tsukamoto, K., & Uchida, K. (1992). Migration mechanism of the ayu. In Oceanic and Anthropogenic Controls of Life in the Pacific Ocean: Proceedings of the 2nd Pacific Symposium on Marine Sciences, Nadhodka, Russia, August 11–19, 1988 (pp. 145-172). Dordrecht: Springer Netherlands. https://doi.org/10.1007/978-94-011-2773-8_12
[82] Yoshioka, H. (2016). Mathematical analysis and validation of an exactly solvable model for upstream migration of fish schools in one-dimensional rivers. Mathematical biosciences, 281, 139-148. https://doi.org/10.1016/j.mbs.2016.09.014
[83] Bjerck, H. B., Urke, H. A., Haugen, T. O., Alfredsen, J. A., Ulvund, J. B., & Kristensen, T. (2021). Synchrony and multimodality in the timing of Atlantic salmon smolt migration in two Norwegian fjords. Scientific Reports, 11(1), 6504. https://doi.org/10.1038/s41598-021-85941-9
[84] Kao, A. B., Banerjee, S. C., Francisco, F. A., & Berdahl, A. M. (2024). Timing decisions as the next frontier for collective intelligence. Trends in Ecology & Evolution, 39(10), 904-912. https://doi.org/10.1016/j.tree.2024.06.003
[85] Van Wichelen, J., Buysse, D., Verhelst, P., Belpaire, C., Goegebeur, M., Vlietinck, K., & Coeck, J. (2023). Nocturnal tidal barrier management improves glass eel migration in times of drought and salinization risk. River Research and Applications, 39(4), 797-801. https://doi.org/10.1002/rra.4088
[86] Verhelst, P., Westerberg, H., Coeck, J., Harrison, L., Moens, T., Reubens, J., ... & Righton, D. (2023). Tidal and circadian patterns of European eel during their spawning migration in the North Sea and the English Channel. Science of the total environment, 905, 167341. https://doi.org/10.1016/j.scitotenv.2023.167341
[87] Mei, Y., Al-Jarrah, M., Taghvaei, A., & Chen, Y. (2025). Flow matching for stochastic linear control systems. Proceedings of the 7th Annual Learning for Dynamics & Control Conference, PMLR 283:484-496, 2025. https://proceedings.mlr.press/v283/mei25a.html
[88] Caradima, B., Scheidegger, A., Brodersen, J., & Schuwirth, N. (2021). Bridging mechanistic conceptual models and statistical species distribution models of riverine fish. Ecological Modelling, 457, 109680. https://doi.org/10.1016/j.ecolmodel.2021.109680
[89] Filar, J. A., Holden, M. H., Mendiolar, M., & Streipert, S. H. (2025). Persistence index for harvested populations. Mathematical Biosciences, 109497. https://doi.org/10.1016/j.mbs.2025.109497
[90] Nugroho, S. (2025). Operationalizing social-ecological system-based fishery management employing a system dynamics model: Lessons from eel fishery. Ecological Modelling, 509, 111276. https://doi.org/10.1016/j.ecolmodel.2025.111276
[91] Karatzas, I., & Shreve, S. (1998). Brownian motion and stochastic calculus. Springer, New York.
[92] Alfonsi, A. (2015). Affine diffusions and related processes: simulation, theory and applications. Springer, Cham.
[93] Mao, X. (2007). Stochastic differential equations and applications. Elsevier.

[94] Okayama Prefectural Fisheries Experimental Station (2008). Survey of wild sweetfish migration https://www.pref.okayama.jp/uploaded/attachment/350487.pdf (in Japanese, last accessed on August 16, 2025)
[95] Okayama Prefectural Agriculture, Forestry and Fisheries Center (2014). https://www.pref.okayama.jp/uploaded/attachment/206364.pdf (in Japanese, last accessed on August 16, 2025)
[96] Yamagata Prefecture (2014). A simple method for estimating the number of sweetfish migrating upstream through the Nagasawa Dam fishway on the Oguni River https://www.pref.yamagata.jp/documents/6277/h27-yamagata-naisuishi-seika04.pdf (in Japanese, last accessed on August 16, 2025)
[97] Azuma, K. et al. (2020) Counting the numbers of ascending juvenile ayu based on visual observations and underwater videos in the Shimanto River, western Japan. Aquaculture Science, 68(4), p.375-382. In Japanese with English Abstract. https://doi.org/10.11233/aquaculturesci.68.375
[98] Takahashi, N., Kitamura, Y., & Shimizu, K. (2012). Evaluation of the Migration Environment of Ayu under the Different Flow Conditions — A case of the Eino weir in the Hatto River, Sendai River system —. Transactions of The Japanese Society of Irrigation, Drainage and Rural Engineering, 80(4), 305-310. In Japanese with English Abstract. https://doi.org/10.11408/jsidre.80.305
[99] Shimada, K., Goto, K., Yamamoto, K., & Wada, Y. (2006). Forecasting ascending population sizes of amphidromous juvenile ayu *Plecoglossus altivelis* in the Nagara River. Nippon Suisan Gakkaishi, 72(4), 665-672. In Japanese with English Abstract. https://doi.org/10.2331/suisan.72.665